\documentclass[11pt]{article}
\usepackage{amsmath}
\usepackage{amssymb}
\newcommand\CC{\mathbb{C}}
\newcommand\RR{\mathbb{R}}
\newcommand\ZZ{\mathbb{Z}}
\newcommand\al\alpha
\newcommand\be\beta
\newcommand\ga\gamma
\newcommand\de\delta
\newcommand\tha\theta
\newcommand\la\lambda
\newcommand\si\sigma
\newcommand\om\omega
\newcommand\Ga{\Gamma}
\newcommand\De{\Delta}
\newcommand\half{\frac12}
\newcommand\thalf{\tfrac12}
\newcommand\iy\infty
\newcommand\wt\widetilde
\newcommand\Zpos{\ZZ_{>0}}
\newcommand\Znonneg{\ZZ_{\ge0}}
\newcommand\lan{\langle}
\newcommand\ran{\rangle}
\renewcommand\Re{{\rm Re}\,}
\renewcommand\Im{{\rm Im}\,}
\newcommand{\hyp}[5]{\,\mbox{}_{#1}F_{#2}\left(
  \genfrac{}{}{0pt}{}{#3}{#4};#5\right)}
\newcommand\LHS{left-hand side}
\newcommand\RHS{right-hand side}
\newcommand{\qhyp}[5]{\,\mbox{}_{#1}\phi_{#2}\left[
  \genfrac{}{}{0pt}{}{#3}{#4};#5\right]}
\newcommand{\qbihyp}[5]{\,\mbox{}_{#1}\psi_{#2}\left[
  \genfrac{}{}{0pt}{}{#3}{#4};#5\right]}
\newcommand\qbinom[3]{{\genfrac[]{0pt}0{#1}{#2}}_#3}
\newcommand\td{\tilde}
\newcommand{\dup}{\textup{d}}
\newcommand{\eup}{\textup{e}}
\newcommand{\iup}{\mkern1mu\textup{i}\mkern1mu}

\usepackage[pstarrows]{pict2e}

\numberwithin{equation}{section}
\begin{document}
\title{$q$-Special functions, an overview}
\author{Tom H. Koornwinder\footnote{Korteweg-de Vries Institute,
University of Amsterdam, P.O.\ Box 94248, 1090 GE Amsterdam, The Netherlands;
email: \texttt{thkmath@xs4all.nl}}}
\date{}
\maketitle
\begin{abstract}
This article gives a brief introduction to $q$-special functions, i.e.,
$q$-analogues of the classical special functions. Here $q$ is a deformation
parameter, usually $0<q<1$, where $q=1$ is the classical case. The main topics
to be treated are $q$-hypergeometric series, with some selected evaluation and
transformation formulas, and the $q$-hypergeometric orthogonal polynomials,
most notably the Askey--Wilson polynomials. Some newer topics as nonsymmetric
analogues and $q=-1$ limits will also be addressed. In several variables we
discuss Macdonald polynomials associated with root systems, in particular the
$A_n$ and the $BC_n$ case. The theory of elliptic hypergeometric series also
gets some attention. The occurrence of $q$-series in number theory and
combinatorics will be discussed. Finally we indicate applications and
interpretations in quantum groups, Chevalley groups, affine Lie algebras and
statistical mechanics.
\end{abstract}
\newpage
\section{Introduction}
\label{sec:06}
Special functions originated as solutions of the
classical PDEs of mathematical physics by separation of variables.
Often they formed orthogonal systems
such as Legendre polynomials,
which were building blocks for obtaining more general solutions.
Soon it was observed that such functions
fitted into function classes having many parameters: the
hypergeometric series.
An additional parameter $q$ first occurred in Euler's generating functions for
partition numbers (see \S\ref{122}). This parameter $q$ could also be considered
as a deformation parameter giving the classical case for $q=1$. In this way,
the fact that $(1-q^a)(1-q^{a+1})\ldots(1-q^{a+k-1})/(1-q)^k$ is a deformation of
the shifted factorial $a(a+1)\ldots(a+k-1)$ gave rise to the definitions
of the $q$-shifted factorial \eqref{eq:01.01} and the $q$-hypergeometric series
\eqref{eq:01.03}. Parallel to this was a $q$-deformation of the calculus,
leading to the $q$-derivative and the $q$-integral, see \S\ref{123}.
Surprising summation and transformation formulas were found for special
$q$-hypergeometric series, see \S\ref{124}. Furthermore, some orthogonal polynomials
expressible as $q$-hypergeometric series were found.

The field got an enormous boost from 1975 onward by the efforts of Askey and
coworkers. This culminated into the introduction of the Askey--Wilson
polynomials \cite{29} and the publication of the monographs \cite{14} and \cite{01}.
The root systems from semisimple Lie theory were a fruitful guide
for Macdonald to obtain analogues in several variables of the $q$-hypergeometric
orthogonal polynomials, see Section \ref{sec:05}. These polynomials are symmetric
under the Weyl group for the corresponding root system.
Dunkl's differential-reflection
operator got a $q$-analogue in the framework of Cherednik's double affine
Hecke algebra, thus leading to non-symmetric polynomials and to the solution
of Macdonald's problems for the symmetric polynomials. Elliptic
hypergeometric series were introduced as a next level above the $q$-level.
The deformation parameter $q$, in use for more than two centuries, got
a new meaning as $q$ for quantum by the introduction of quantum groups
in the period 1980--1990. Some $q$-hypergeometric orthogonal polynomials and
their analogues in several variables turned out to live on quantum groups.

\paragraph{Conventions}\quad\\
$q\in\CC\backslash\{1\}$ in general, but $|q|<1$ in all infinite sums and products
and $0<q<1$ for $q$-hypergeometric orthogonal polynomials.\\
$n,m,N$ will be nonnegative integers unless mentioned otherwise.
\newpage
\tableofcontents
\section{$q$-Hypergeometric series}
\label{sec:01}
\subsection{Definitions}
Standard reference for Section \ref{sec:01} is Gasper \& Rahman \cite{01}.\\
For $a\in\CC$ the {\em $q$-shifted factorial} $(a;q)_k$ is defined as a product of
$k$ factors:
\begin{equation}
(a;q)_k:=(1-a)(1-aq)\ldots(1-aq^{k-1})\quad(k\in\Zpos);\qquad
(a;q)_0:=1.
\label{eq:01.01}
\end{equation}
If $|q|<1$ this definition remains meaningful for $k=\iy$ as a convergent
infinite product:
\begin{equation}
(a;q)_\iy:=\prod_{j=0}^\iy(1-aq^j).
\label{eq:01.02}
\end{equation}
We also write $(a_1,\ldots,a_r;q)_k$ for the product of $r$
$q$-shifted factorials:
\begin{equation}
(a_1,\ldots,a_r;q)_k:=(a_1;q)_k\ldots(a_r;q)_k\quad
\mbox{($k\in\Znonneg$ or $k=\iy$).}
\label{eq:01.03}
\end{equation}

A {\em $q$-hypergeometric series} is a power series
(for the moment still formal)
in one complex variable
$z$ with power series coefficients which depend, apart from $q$, on
$r$ complex {\em upper parameters} $a_1,\ldots,a_r$ and $s$ complex
{\em lower parameters} $b_1,\ldots,b_s$ as follows:
\begin{multline}
\qhyp rs{a_1,\ldots,a_r}{b_1,\ldots,b_s}{q,z}=
{}_r\phi_s(a_1,\ldots,a_r;b_1,\ldots,b_s;q,z)\\
:=\sum_{k=0}^\iy\frac{(a_1,\ldots,a_r;q)_k}{(b_1,\ldots,b_s;q)_k\,(q;q)_k}\,
\left((-1)^k q^{\half k(k-1)}\right)^{s-r+1} z^k\qquad(r,s\in\Znonneg).
\label{eq:01.04}
\end{multline}
Clearly the above expression is symmetric in $a_1,\ldots,a_r$ and symmetric
in $b_1,\ldots,b_s$.
On the \RHS\ of \eqref{eq:01.04} we have that
\begin{equation}
\frac{\mbox{($k+1$)th term}}{\mbox{$k$th term}}=
\frac{(1-a_1q^k)\ldots(1-a_rq^k)\,(-q^k)^{s-r+1}\,z}
{(1-b_1q^k)\ldots(1-b_sq^k)\,(1-q^{k+1})}
\label{eq:01.05}
\end{equation}
is rational in $q^k$.
Conversely, any rational function in $q^k$ can be written in the form of
the \RHS\ of \eqref{eq:01.05}.
Hence, any series $\sum_{k=0}^\iy c_k$ with $c_0=1$ and $c_{k+1}/c_k$
rational in $q^k$ is of the form of a $q$-hypergeometric series
\eqref{eq:01.04}.

In order to avoid singularities in the terms of \eqref{eq:01.04} we assume
that $b_1,\ldots,b_s\ne 1,q^{-1},q^{-2},\ldots\,$. If, for some $i$,
$a_i=q^{-n}$ then all terms in the series
\eqref{eq:01.04} with $k>n$ will vanish.
If none of the $a_i$ is equal to $q^{-n}$ and if
$|q|<1$ then
the radius of convergence of the power series \eqref{eq:01.04} equals
$\iy$ if $r<s+1$, $1$ if $r=s+1$, and $0$ if $r>s+1$.

We can view the $q$-shifted factorial as a $q$-analogue of the
{\em shifted factorial} (or {\em Pochhammer symbol})
by the limit formula
\begin{equation}
\lim_{q\to1}\frac{(q^a;q)_k}{(1-q)^k}=(a)_k:=a(a+1)\ldots(a+k-1).
\label{eq:01.06}
\end{equation}
Hence the {\em $q$-binomial coefficient}
\begin{equation}
\qbinom nkq:=\frac{(q;q)_n}{(q;q)_k(q;q)_{n-k}}\qquad(n,k\in\ZZ,\;n\ge k\ge0)
\label{eq:01.07}
\end{equation}
tends to the binomial coefficient for $q\to 1$:
\begin{equation}
\lim_{q\to1}\qbinom nkq=\binom nk,
\label{eq:01.08}
\end{equation}
and a suitably renormalized $q$-hypergeometric series tends (at least formally)
to a hypergeometric series as $q\uparrow1$:
\begin{multline}
\lim_{q\uparrow1}
\qhyp{r+r'}{s+s'}{q^{a_1},\ldots,q^{a_r},c_1,\ldots,c_{r'}}
{q^{b_1},\ldots,q^{b_s},d_1,\ldots,d_{s'}}{q,(q-1)^{1+s-r}z}\\
=\hyp rs{a_1,\ldots,a_r}{b_1,\ldots,b_s}
{\frac{(c_1-1)\ldots(c_{r'}-1)\,z}{(d_1-1)\ldots(d_{s'}-1)}}.
\label{eq:01.09}
\end{multline}

At least formally, there are limit relations between $q$-hypergeometric
series with neighbouring $r,s$:
\begin{align}
\lim_{a_r\to\iy}\qhyp rs{a_1,\ldots,a_r}{b_1,\ldots,b_s}{q,\frac z{a_r}}
&=\qhyp{r-1}s{a_1,\ldots,a_{r-1}}{b_1,\ldots,b_s}{q,z},
\label{eq:01.17}\\
\lim_{b_s\to\iy}\qhyp rs{a_1,\ldots,a_r}{b_1,\ldots,b_s}{q,b_sz}
&=\qhyp r{s-1}{a_1,\ldots,a_r}{b_1,\ldots,b_{s-1}}{q,z}.
\label{eq:01.18}
\end{align}

A terminating $q$-hypergeometric series $\sum_{k=0}^n c_k\,z^k$
rewritten as
$z^n\sum_{k=0}^n c_{n-k} z^{-k}$ yields another terminating
$q$-hypergeometric series, for instance:
\begin{multline}
\qhyp{s+1}s{q^{-n},a_1,\ldots,a_s}{b_1,\ldots,b_s}{q,z}
=(-1)^n\,q^{-\half n(n+1)}\,
\frac{(a_1,\ldots,a_n;q)_n}{(b_1,\ldots,b_s;q)_n}\,z^n\\
\times\qhyp{s+1}s{q^{-n},q^{-n+1}b_1^{-1},\ldots,q^{-n+1}b_s^{-1}}
{q^{-n+1}a_1^{-1},\ldots,q^{-n+1}a_s^{-1}}
{q,\frac{q^{n+1}b_1\ldots b_s}{a_1\ldots a_s z}}.
\label{eq:01.60}
\end{multline}

Often, in physics and quantum groups related literature, the following
notation is used for $q$-number, $q$-factorial and $q$-Pochhammer symbol:
\begin{equation}
[a]_q:=\frac{q^{\half a}-q^{-\half a}}{q^\half-q^{-\half}},\quad
[k]_q!:=\prod_{j=1}^k[j]_q,\quad
([a]_q)_k:=\prod_{j=0}^{k-1}[a+j]_q\qquad(k\in\Znonneg).
\label{eq:01.10}
\end{equation}
For $q\to1$ these symbols tend to their classical counterparts without
the need for renormalization. They are expressed in terms of the
standard notation \eqref{eq:01.01} as follows:
\begin{equation}
[k]_q!=q^{-\frac14 k(k-1)}\,\frac{(q;q)_k}{(1-q)^k},\qquad
([a]_q)_k=q^{-\half k(a-1)}\,q^{-\frac14 k(k-1)}\,
\frac{(q^a;q)_k}{(1-q)^k}\,.
\label{eq:01.11}
\end{equation}
\subsection{Special cases}
For $s=r-1$ formula \eqref{eq:01.04} simplifies to
\begin{equation}
\qhyp r{r-1}{a_1,\ldots,a_r}{b_1,\ldots,b_{r-1}}{q,z}=
\sum_{k=0}^\iy\frac{(a_1,\ldots,a_r;q)_k}{(b_1,\ldots,b_{r-1};q)_k\,(q;q)_k}\,
z^k,
\label{eq:01.12}
\end{equation}
which has radius of convergence 1 in the non-terminating case.
The case $r=2$ of \eqref{eq:01.12} is the $q$-analogue of the Gauss
hypergeometric series.
\paragraph{$q$-Binomial series}
\begin{equation}
{}_1\phi_0(a;-;q,z)=
\sum_{k=0}^\iy\frac{(a;q)_k z^k}{(q;q)_k}=
\frac{(az;q)_\iy}{(z;q)_\iy}\quad
\mbox{(if series is not terminating then $|z|<1$).}
\label{eq:01.13}
\end{equation}
\paragraph{$q$-Exponential series}
\begin{align}
e_q(z):=&{}_1\phi_0(0;-;q,z)=\sum_{k=0}^\iy\frac{z^k}{(q;q)_k}=
\frac1{(z;q)_\iy}\qquad(|z|<1),
\label{eq:01.14}\\
E_q(z):=&{}_0\phi_0(-;-;q,-z)=\sum_{k=0}^\iy\frac{q^{\half k(k-1)}z^k}{(q;q)_k}
=(-z;q)_\iy=\left(e_q(-z)\right)^{-1}\qquad(z\in\CC),
\label{eq:01.15}\\
\varepsilon_q(z):=&{}_1\phi_1(0;-q^\half;q^\half,-z)=
\sum_{k=0}^\iy \frac{q^{\frac14 k(k-1)}}{(q;q)_k}\,z^k\qquad(z\in\CC;\;
\mbox{notation not standard}).
\label{eq:01.64}
\end{align}
\paragraph{Jackson's $q$-Bessel functions}
\begin{align}
J_\nu^{(1)}(x;q):=&\frac{(q^{\nu+1};q)_\iy}{(q;q)_\iy}\,(\thalf x)^\nu\,
\qhyp21{0,0}{q^{\nu+1}}{q,-\tfrac14 x^2}\qquad(0<x<2),
\label{eq:01.61}\\
J_\nu^{(2)}(x;q):=&\frac{(q^{\nu+1};q)_\iy}{(q;q)_\iy}\,(\thalf x)^\nu\,
\qhyp01{-}{q^{\nu+1}}{q,-\tfrac14 q^{\nu+1}x^2}=
(-\tfrac14x;q)_\iy\,J_\nu^{(1)}(x;q)\quad(x>0),
\label{eq:01.62}\\
J_\nu^{(3)}(x;q):=&\frac{(q^{\nu+1};q)_\iy}{(q;q)_\iy}\,(\thalf x)^\nu\,
\qhyp11{0}{q^{\nu+1}}{q,\tfrac14 q x^2}\qquad(x>0).
\label{eq:01.63}
\end{align}
See \eqref{eq:02.35} for the orthogonality relation for $J_\nu^{(3)}(x;q)$.

If $\exp_q(z)$ denotes one of the three $q$-exponentials
\eqref{eq:01.14}--\eqref{eq:01.64} then
$\thalf\bigl(\exp_q(ix)+\exp_q(-ix)\bigr)$ is a $q$-analogue of the cosine and
$-\thalf i\bigl(\exp_q(ix)-\exp_q(-ix)\bigr)$ is a $q$-analogue of the sine.
The three $q$-cosines are essentially the case $\nu=-\thalf$ of
the corresponding $q$-Bessel functions \eqref{eq:01.61}--\eqref{eq:01.63},
and the three $q$-sines are essentially the case $\nu=\thalf$ of
$x$ times the corresponding $q$-Bessel functions.
See also Suslov \cite{05}.
\subsection{$q$-Derivative and $q$-integral}
\label{123}
The {\em $q$-derivative} of a function $f$ given on a subset of $\RR$ or
$\CC$ is defined by
\begin{equation}
(D_qf)(x):=\frac{f(x)-f(qx)}{(1-q)x}\qquad(x\ne0,\;q\ne1),
\label{eq:01.19}
\end{equation}
where $x$ and $qx$ should be in the domain of $f$. By continuity we set
$(D_qf)(0):=f'(0)$, provided $f'(0)$ exisits. If $f$ is differentiable
on an open interval $I$ then
\begin{equation}
\lim_{q\uparrow1}(D_qf)(x)=f'(x)\qquad(x\in I).
\label{eq:01.20}
\end{equation}

For $a\in\RR\backslash\{0\}$ and a function $f$ given on $(0,a]$
or $[a,0)$, we define the {\em $q$-integral} by
\begin{equation}
\int_0^a f(x)\,\dup_qx:=a(1-q)\sum_{k=0}^\iy f(aq^k)\,q^k=
\sum_{k=0}^\iy f(aq^k)\,(aq^k-aq^{k+1}),
\label{eq:01.21}
\end{equation}
provided the infinite sum converges absolutely (for instance if $f$ is
bounded). If $F(a)$ is given by the \LHS\ of \eqref{eq:01.21} then
$D_qF=f$. The \RHS\ of \eqref{eq:01.21} is an infinite Riemann sum.
For $q\uparrow 1$ it converges, at least formally, to
$\int_0^a f(x)\,\dup x$.

For nonzero $a,b\in\RR$ we define
\begin{equation}
\int_a^b f(x)\,\dup_qx:=\int_0^b f(x)\,\dup_qx-\int_0^a f(x)\,\dup_qx.
\label{eq:01.22}
\end{equation}
For a $q$-integral over $(0,\iy)$ we have to specify a $q$-lattice
$\{aq^k\}_{k\in\ZZ}$ for some $a>0$ (up to multiplication by an integer
power of $q$):
\begin{equation}
\int_0^{a.\iy}f(x)\,\dup_qx:=a(1-q)\sum_{k=-\iy}^\iy f(aq^k)\,q^k=
\lim_{n\to\iy}\int_0^{q^{-n}a}f(x)\,\dup_qx.
\label{eq:01.23}
\end{equation}
\subsection{The $q$-gamma and $q$-beta functions}
The {\em $q$-gamma function} is defined by
\begin{align}
\Ga_q(z):=&\frac{(q;q)_\iy\,(1-q)^{1-z}}{(q^z;q)_\iy}\qquad(z\ne0,-1,-2,\ldots)
\label{eq:01.24}\\
=&\int_0^{(1-q)^{-1}}t^{z-1}\,E_q(-(1-q)qt)\,\dup_qt\qquad(\Re z>0)
\label{eq:01.25}.
\end{align}
Then
\begin{align}
\Ga_q(z+1)&=\frac{1-q^z}{1-q}\,\Ga_q(z),
\label{eq:01.26}\\
\Ga_q(n+1)&=\frac{(q;q)_n}{(1-q)^n}\,,
\label{eq:01.27}\\
\lim_{q\uparrow1}\Ga_q(z)&=\Ga(z).
\label{eq:01.28}
\end{align}

The {\em $q$-beta function} is defined by
\begin{align}
B_q(a,b):=&\frac{\Ga_q(a)\Ga_q(b)}{\Ga_q(a+b)}
=\frac{(1-q)\,(q,q^{a+b};q)_\iy}{(q^a,q^b;q)_\iy}\qquad
(a,b\ne0,-1,-2,\ldots),
\label{eq:01.29}\\
=&\int_0^1 t^{b-1}\,\frac{(qt;q)_\iy}{(q^a t;q)_\iy}\,\dup_qt\qquad
(\Re b>0,\;a\ne0,-1,-2,\ldots).
\label{eq:01.30}
\end{align}
\subsection{The $q$-Gauss hypergeometric series}
\paragraph{$q$-Analogue of Euler's integral representation}
\begin{equation}
{}_2\phi_1(q^a,q^b;q^c;q,z)=
\frac{\Ga_q(c)}{\Ga_q(a)\Ga_q(c-b)}\,
\int_0^1 t^{b-1}\,\frac{(tq;q)_\iy}{(tq^{c-b};q)_\iy}\,
\frac{(tzq^a;q)_\iy}{tz;q)_\iy}\,\dup_qt\quad(\Re b>0,\;|z|<1).
\label{eq:01.31}
\end{equation}
By substitution of \eqref{eq:01.21}, formula \eqref{eq:01.31} becomes
a transformation formula:
\begin{equation}
{}_2\phi_1(a,b;c;q,z)=\frac{(az;q)_\iy}{(z;q)_\iy}\,
\frac{(b;q)_\iy}{(c;q)_\iy}\,{}_2\phi_1(c/b,z;az;q,b).
\label{eq:01.35}
\end{equation}
Note the mixing of argument $z$ and parameters $a,b,c$ on the \RHS.
\paragraph{Evaluation formulas in special points}
\begin{align}
{}_2\phi_1\bigl(a,b;c;q,c/(ab)\bigr)&=
\frac{(c/a,c/b;q)_\iy}{(c,c/(ab);q)_\iy}\qquad
(|c/(ab)|<1),
\label{eq:01.32}\\
{}_2\phi_1(q^{-n},b;c;q,cq^n/b)&=\frac{(c/b;q)_n}{(c;q)_n}\,,
\label{eq:01.33}\\
{}_2\phi_1(q^{-n},b;c;q,q)&=\frac{(c/b;q)_n\,b^n}{(c;q)_n}\,.
\label{eq:01.34}
\end{align}
\paragraph{Two general transformation formulas}
\begin{align}
\qhyp21{a,b}c{q,z}&=
\frac{(az;q)_\iy}{(z;q)_\iy}\,
\qhyp22{a,c/b}{c,az}{q,bz},
\label{eq:01.36}\\
&=\frac{(abz/c;q)_\iy}{(z;q)_\iy}\,
\qhyp21{c/a,c/b}c{q,\frac{abz}c}.
\label{eq:01.37}
\end{align}
\paragraph{Transformation formulas in the terminating case}
\begin{align}
\qhyp21{q^{-n},b}c{q,z}
&=\frac{(c/b;q)_n}{(c;q)_n}\,
\qhyp32{q^{-n},b,q^{-n}bc^{-1}z}{q^{1-n}bc^{-1},0}{q,q}
\label{eq:01.38}\\
&=(q^{-n}bc^{-1}z;q)_n\,
\qhyp32{q^{-n},cb^{-1},0}{c,qcb^{-1}z^{-1}}{q,q}
\label{eq:01.39}\\
&=\frac{(c/b;q)_n}{(c;q)_n}\,b^n\,
\qhyp31{q^{-n},b,qz^{-1}}{q^{1-n}bc^{-1}}{q,\frac zc}.
\label{eq:01.40}
\end{align}

\paragraph{Second order $q$-difference equation}
\begin{multline}
z(q^c-q^{a+b+1}z)(D_q^2u)(z)+
\left(\frac{1-q^c}{1-q}-
\left(q^b\frac{1-q^a}{1-q}+q^a\frac{1-q^{b+1}}{1-q}\right)z\right)
(D_qu)(z)\\
-\frac{1-q^a}{1-q}\,\frac{1-q^b}{1-q}\,u(z)=0.
\label{eq:01.41}
\end{multline}
Some special solutions of \eqref{eq:01.41} are:
\begin{align}
u_1(z)&:={}_2\phi_1(q^a,q^b;q^c;q,z),
\label{eq:01.42}\\
u_2(z)&:=z^{1-c}\,{}_2\phi_1(q^{1+a-c},q^{1+b-c};q^{2-c};q,z),
\label{eq:01.43}\\
u_3(z)&:=z^{-a}\,{}_2\phi_1(q^a,q^{a-c+1};q^{a-b+1};q,q^{-a-b+c+1}z^{-1}).
\label{eq:01.44}
\end{align}
They are related by:
\begin{multline}
u_1(z)+
\frac{(q^a,q^{1-c},q^{c-b};q)_\iy}{(q^{c-1},q^{a-c+1},q^{1-b};q)_\iy}\,
\frac{(q^{b-1}z,q^{2-b}z^{-1};q)_\iy}{(q^{b-c}z,q^{c-b+1}z^{-1};q)_\iy}\,
u_2(z)\\
=\frac{(q^{1-c},q^{a-b+1};q)_\iy}{(q^{1-b},q^{a-c+1};q)_\iy}\,
\frac{(q^{a+b-c}z,q^{c-a-b+1}z^{-1};q)_\iy\,z^a}
{(q^{b-c}z,q^{c-b+1}z^{-1};q)_\iy}\,u_3(z).
\label{eq:01.45}
\end{multline}
\subsection{Summation and transformation formulas for ${}_r\phi_{r-1}$ series}
\label{124}
An ${}_r\phi_{r-1}$ series \eqref{eq:01.12} is called {\em balanced}
if $b_1\ldots b_{r-1}=q a_1\ldots a_r$ and $z=q$, and the series
is called {\em very-well-poised}
if $qa_1=a_2b_1=a_3b_2=\cdots=a_rb_{r-1}$ and $qa_1^\half=a_2=-a_3$.
The following more compact notation is used for very-well-poised series:
\begin{equation}
{}_rW_{r-1}(a_1;a_4,a_5,\ldots,a_r;q,z):=
\qhyp r{r-1}{a_1,qa_1^\half,-qa_1^\half,a_4,\ldots,a_r}
{a_1^\half,-a_1^\half,qa_1/a_4,\ldots,qa_1/a_r}{q,z}.
\label{eq:01.57}
\end{equation}
Below only a few of the most important identities are given.
See \cite{01} for many more. An important tool for obtaining
complicated identities from more simple ones is {\em Bailey's Lemma},
which can moreover be iterated ({\em Bailey chain}), see
\cite[Ch.3]{19}.
\paragraph{The $q$-Saalsch\"utz sum for a terminating balanced ${}_3\phi_2$}
\begin{equation}
\qhyp32{a,b,q^{-n}}{c,q^{1-n}abc^{-1}}{q,q}=
\frac{(c/a,c/b;q)_n}{(c,c/(ab);q)_n}\,.
\label{eq:01.46}
\end{equation}
\paragraph{Jackson's sum for a terminating
balanced ${}_8W_7$}
\begin{equation}
{}_8W_7(a;b,c,d,q^{n+1}a^2/(bcd),q^{-n};q,q)
=\frac{(qa,qa/(bc),qa/(bd),qa/(cd);q)_n}{(qa/b,qa/c,qa/d,qa/(bcd);q)_n}\,.
\label{eq:01.47}
\end{equation}
\paragraph{Watson's transformation of a
terminating ${}_8W_7$ into a terminating balanced ${}_4\phi_3$}
\begin{equation}
{}_8W_7\left(a;b,c,d,e,q^{-n};q,\frac{q^{n+2}a^2}{bcde}\right)
=\frac{(qa,qa/(de);q)_n}{(qa/d,qa/e;q)_n}\,
\qhyp43{q^{-n},d,e,qa/(bc)}{qa/b,qa/c,q^{-n}de/a}{q,q}.
\label{eq:01.48}
\end{equation}
\paragraph{Sears' transformation of a terminating balanced ${}_4\phi_3$}\begin{equation}
\qhyp43{q^{-n},a,b,c}{d,e,f}{q,q}=
\frac{(e/a,f/a;q)_n}{(e,f;q)_n}\,a^n\,
\qhyp43{q^{-n},a,d/b,d/c}{d,q^{1-n}a/e,q^{1-n}a/f}{q,q}.
\label{eq:01.65}
\end{equation}
By iteration and by symmetries in the upper and in the lower parameters,
many other versions of this identity can be found. An elegant comprehensive
formulation of all these versions is as follows.\\
Let $x_1x_2x_3x_4x_5x_6=q^{1-n}$.
Then the following expression is symmetric in
$x_1,x_2,x_3,x_4,x_5,x_6$:
\begin{equation}
\frac{q^{\half n(n-1)}(x_1x_2x_3x_4,x_1x_2x_3x_5,x_1x_2x_3x_6;q)_n}
{(x_1x_2x_3)^n}\,
\qhyp43{q^{-n},x_2x_3,x_1x_3,x_1x_2}{x_1x_2x_3x_4,x_1x_2x_3x_5,x_1x_2x_3x_6}
{q,q}.
\label{eq:01.59}
\end{equation}
Similar formulations involving symmetry groups can be given for other
transformations, see \cite{03}.
\paragraph{Bailey's transformation of a terminating balanced ${}_{10}W_9$}
\begin{multline}
{}_{10}W_9\left(a;b,c,d,e,f,\frac{q^{n+2}a^3}{bcdef},q^{-n};q,q\right)
=\frac{(qa,qa/(ef),(qa)^2/(bcde),(qa)^2/(bcdf);q)_n}
{(qa/e,qa/f,(qa)^2/(bcdef),(qa)^2/(bcd);q)_n}\\
\times
{}_{10}W_9\left(\frac{qa^2}{bcd};\frac{qa}{cd},\frac{qa}{bd},\frac{qa}{bc},
e,f,\frac{q^{n+2}a^3}{bcdef},q^{-n};q,q\right).
\label{eq:01.58}
\end{multline}
\subsection{Rogers--Ramanujan identities}
\begin{align}
{}_0\phi_1(-;0;q,q)=\sum_{k=0}^\iy\frac{q^{k^2}}{(q;q)_k}=&
\frac1{(q,q^4;q^5)_\iy}\,,
\label{eq:01.49}\\
{}_0\phi_1(-;0;q,q^2)=\sum_{k=0}^\iy\frac{q^{k(k+1)}}{(q;q)_k}=&
\frac1{(q^2,q^3;q^5)_\iy}\,.
\label{eq:01.50}
\end{align}
\subsection{Bilateral series}
Definition \eqref{eq:01.01} can be extended by
\begin{equation}
(a;q)_k:=\frac{(a;q)_\iy}{(aq^k;q)_\iy}\qquad(k\in\ZZ).
\label{eq:01.51}
\end{equation}
Define a {\em bilateral $q$-hypergeometric series} by the Laurent series
\begin{multline}
\qbihyp rs{a_1,\ldots,a_r}{b_1,\ldots,b_s}{q,z}
={}_r\psi_s(a_1,\ldots,a_r;b_1,\ldots,b_s;q,z)\\
:=
\sum_{k=-\iy}^\iy\frac{(a_1,\ldots,a_r;q)_k}{(b_1,\ldots,b_s;q)_k}
\left((-1)^k q^{\half k(k-1)}\right)^{s-r} z^k\quad
(a_1,\ldots,a_r,b_1,\ldots,b_s\ne0,\;s\ge r).
\label{eq:01.52}
\end{multline}
The Laurent series is convergent if
$|b_1\ldots b_s/(a_1\ldots a_r)|\allowbreak<|z|$
and moreover, for $s=r$, $|z|<1$.
\paragraph{Ramanajan's ${}_1\psi_1$ summation formula}
\begin{equation}
{}_1\psi_1(b;c;q,z)=\frac{(q,c/b,bz,q/(bz);q)_\iy}
{(c,q/b,z,c/(bz);q)_\iy}\qquad(|c/b|<|z|<1).
\label{eq:01.53}
\end{equation}
This has as a limit case
\begin{equation}
{}_0\psi_1(-;c;q,z)=\frac{(q,z,q/z;q)_\iy}{(c,c/z;q)_\iy}\quad
(|z|>|c|),
\label{eq:01.54}
\end{equation}
and as a further specialization the
{\em Jacobi triple product identity}
\begin{equation}
\sum_{k=-\iy}^\iy(-1)^k\,q^{\half k(k-1)}\,z^k=(q,z,q/z;q)_\iy\qquad
(z\ne0),
\label{eq:01.55}
\end{equation}
which can be rewritten as a product formula for a {\em theta function}:
\begin{equation}
\tha_4(x;q):=\sum_{k=-\iy}^\iy(-1)^k\,q^{k^2}\,e^{2\pi ikx}=
\prod_{k=1}^\iy(1-q^{2k})\bigl(1-2q^{k-1}\cos(2\pi x)+q^{4k-2}\bigr).
\label{eq:01.56}
\end{equation}
\section{$q$-Hypergeometric orthogonal polynomials}
\label{sec:02}
In this section, although many formulas will remain valid for more general $q$,
assume $0<q<1$ for orthogonality properties.
We will discuss families of orthogonal polynomials $\{p_n(x)\}$ which
are expressible as terminating $q$-hypergeometric series
and for which either (i) $P_n(x):=p_n(x)$ or (ii)
$P_n(x):=p_n\bigl(\thalf(x+x^{-1})\bigr)$ are eigenfunctions of a
second order $q$-difference operator, i.e.:
\begin{equation}
A(x)\,P_n(qx)+B(x)\,P_n(x)+C(x)\,P_n(q^{-1}x)=\la_n\,P_n(x),
\label{eq:02.01}
\end{equation}
where $A(x)$, $B(x)$ and $C(x)$ are independent of $n$, and where the
$\la_n$ are the eigenvalues.
The generic cases are the four-parameter classes of
{\em Askey--Wilson polynomials}
(continuous weight function) and
{\em $q$-Racah polynomials} (discrete weights
on finitely many points).
They are of type (ii) (quadratic $q$-lattice).
All other cases can be obtained from the generic cases by specialization
or limit transition. In particular, one thus obtains the generic
three-parameter classes of type~(i) (linear $q$-lattice). These are
the {\em big $q$-Jacobi polynomials} (orthogonality by $q$-integral) and
the {\em $q$-Hahn polynomials} (discrete weights
on finitely many points).
For all these families the standard formulas are given in
Koekoek \& Swarttouw \cite[Chapter 14]{02}.
\subsection{$q$-Askey scheme and Verde-Star's description}
Just as the hypergeometric orthogonal polynomials can be arranged in the
Askey scheme \cite[Appendix]{29}, \cite[Chapter 9]{02}, the
$q$-hypergeometric orthogonal polynomials can be arranged in the
\emph{$q$-Askey scheme} \cite[Chapter 14]{02}.
It consists of boxes indicating families of orthogonal
polynomials and arrows between boxes which indicate specializations or
limit transitions. In Figure~\ref{121} the part of the $q$-Askey scheme 
descending from the 4-parameter Askey--Wilson polynomials is drawn. Families
in each next row depend on one parameter less.
\setlength{\unitlength}{2.35mm}
\begin{figure}[ht]
\centering
\hskip-3cm
\begin{picture}(30,40)
\put(15,36)
{\framebox{\begin{tabular}{c}Askey--\\Wilson\end{tabular}}}
\put(16,34) {\vector(-3,-2){14}}
\put(21,34) {\vector(0,-1){2.2}}
\put(23,34) {\vector(5,-1){10}}
\put(15,29)
{\framebox{\begin{tabular}{c}cont.~dual\\$q$-Hahn\end{tabular}}}
\put(21,27) {\vector(0,-1){2.2}}
\put(24,27) {\vector(4,-1){8}}
\put(30,29)
{\framebox{\begin{tabular}{c}big $q$-\\Jacobi\end{tabular}}}
\put(33,27) {\vector(0,-1){2.2}}
\put(36,27) {\vector(4,-1){8}}
\put(-6,22)
{\framebox{\begin{tabular}{c}cont.~$q$-\\Jacobi\end{tabular}}}
\put(-3,20) {\vector(-2,-1){4}}
\put(-1,20) {\vector(2,-1){4}}
\put(17,22)
{\framebox{\begin{tabular}{c}Al-Salam--\\Chihara\end{tabular}}}
\put(21,20) {\vector(-2,-1){4}}
\put(23,20) {\vector(2,-1){4}}
\put(29,22)
{\framebox{\begin{tabular}{c}big $q$-\\Laguerre\end{tabular}}}
\put(32,20) {\vector(-2,-1){4}}
\put(35,20) {\vector(2,-1){4}}
\put(41,22)
{\framebox{\begin{tabular}{c}little $q$-\\Jacobi\end{tabular}}}
\put(45,20) {\vector(-2,-1){4}}
\put(47,20) {\vector(0,-1){3.4}}
\put(-11.5,15)
{\framebox{\begin{tabular}{c}cont.~$q$-\\Laguerre\end{tabular}}}
\put(-4,13) {\vector(3,-1){6}}
\put(-1.5,15)
{\framebox{\begin{tabular}{c}cont.~$q$-ultra-\\spherical\end{tabular}}}
\put(4.7,13) {\vector(0,-1){2.2}}
\put(12,15)
{\framebox{\begin{tabular}{c}cont.~big\\ $q$-Hermite\end{tabular}}}
\put(14.5,13) {\vector(-3,-1){6}}
\put(23,15)
{\framebox{\begin{tabular}{c}Al-Salam--\\Carlitz\end{tabular}}}
\put(28,13) {\vector(0,-1){2.2}}
\put(34.5,15)
{\framebox{\begin{tabular}{c}little $q$-\\Laguerre\end{tabular}}}
\put(37,13) {\vector(3,-1){6}}
\put(45,15)
{\framebox{$q$-Bessel}}
\put(47,14.2) {\vector(-1,-2){1.7}}
\put(1,8)
{\framebox{\begin{tabular}{c}cont.~$q$-\\Hermite\end{tabular}}}
\put(23,8)
{\framebox{\begin{tabular}{c}discrete $q$-\\Hermite\end{tabular}}}
\put(39,8)
{\framebox{\begin{tabular}{c}Stieltjes--\\Wigert\end{tabular}}}
\end{picture}
\vspace*{-1.5cm}
\caption{Part of the $q$-Askey scheme descending from the Askey--Wilson
polynomials}
\label{121}
\end{figure}

Verde-Star \cite{VS2021} (see also \cite{V-Z2016}) described all families in
the $q$-Askey scheme except for the continuous $q$-Hermite polynomials
\eqref{eq:02.13} in a conceptual way as consisting of polynomials
\begin{equation}
u_n(x):=\sum_{k=0}^n\,\prod_{j=0}^{k-1}\frac{(h_n-h_j)(x-x_j)}{g_{j+1}}\,,
\label{110}
\end{equation}
where
\begin{equation}
\begin{split}
&h_k=a_{-1} q^{-k}+a_0+a_1 q^k,\qquad
x_k=b_{-1} q^{-k}+b_0+b_1 q^k,\\
&g_k=d_{-2} q^{-2k}+d_{-1} q^{-k}+d_0+d_1 q^k+d_2 q^{2k},\\
&a_{-1}\notin a_1 q^{\ZZ_{>0}},\quad
\sum_{i=-2}^2 d_i=0,\mbox{ $d_i\ne0$ for some $i$},\\
&d_2=q^{-1}a_1b_1,\quad d_{-2}=qa_{-1}b_{-1}.
\end{split}
\label{109}
\end{equation}
The $q$-Askey scheme can now be redrawn by use of \eqref{109} such that families
are determined by vanishing of one or more extreme terms on the left and on the
right in the expressions for $h_k$, $x_k$ and $g_k$, see \cite{K2022a}.

With explicit data $h_k,x_k,g_k$ for a family in the $q$-Askey scheme
formula \eqref{110} turns down to a $q$-hypergeometric expression.
The formula can also be seen as an explicit expansion of $u_n(x)$ in terms
of \emph{Newton type polynomials}
\begin{equation}
v_k(x):=\prod_{j=0}^{k-1}(x-x_j).
\end{equation}
If an operator $L$ on the space of polynomials is defined by
\begin{equation}
Lv_n:=h_nv_n+g_nv_{n-1}\;(n>0),\quad Lv_0=h_0v_0.
\label{111}
\end{equation}
then the $u_n$ are eigenfunctions of $L$:
\begin{equation}
Lu_n=h_nu_n\quad(n\ge0).
\label{112}
\end{equation}
For each family it turns out that $L$ coincides with a second order
$q$-difference operator.

The above operator $L$ and the operator of multiplication by $x$ both act
on the space of polynomials and generate there an algebra with quite simple
relations for the generators: the \emph{Zhedanov algebra}, see
\cite{Z}, \cite{G-L-Z}. The $q$-Askey scheme can also be redrawn such that families
are determined by vanishing properties of structure coefficients in the
corresponding Zhedanov algebra, see \cite{K2022b}.

\subsection{Askey--Wilson polynomials}
These were introduced by Askey \& Wilson \cite{29},
\paragraph{Definition as $q$-hypergeometric series}
\begin{equation}
p_n(\cos\tha)=p_n(\cos\tha;a,b,c,d\,|\,q):=
\frac{(ab,ac,ad;q)_n}{a^n}\,
R_n(\eup^{\iup\tha};a,b,c,d\,|\,q),
\label{eq:02.02}
\end{equation}
where
\begin{equation}
R_n(z)=R_n(z;a,b,c,d\,|\,q):=
\qhyp43{q^{-n},q^{n-1}abcd,az,az^{-1}}{ab,ac,ad}{q,q}.
\label{101}
\end{equation}
The polynomial $p_n(x;a,b,c,d\,|\,q)$ is symmetric in $a,b,c,d$. For the symmetry
in $a$ and $b$ use Sears' transformation \eqref{eq:01.65}.
\paragraph{Orthogonality relation}
Assume that $a,b,c,d$ are four reals, or two reals and one pair of complex
conjugates, or two pairs of complex conjugates. Also assume that
pairwise products of $a,b,c,d$ are not equal to 1 and have absolute value
$\le1$. Then
\begin{equation}
\int_{-1}^1 p_n(x)\,p_m(x)\,w(x)\,\dup x+\sum_k p_n(x_k)\,p_m(x_k)\,\om_k=
h_n\,\de_{n,m},
\label{eq:02.03}
\end{equation}
where
\begin{equation}
2\pi\sin\tha\,w(\cos\tha)=
\left|\frac{(\eup^{2\iup\tha};q)_\iy}
{(a\eup^{\iup\tha},b\eup^{\iup\tha},c\eup^{\iup\tha},d\eup^{\iup\tha};q)_\iy}\right|^2,
\label{eq:02.04}
\end{equation}
\begin{equation}
h_0=\frac{(abcd;q)_\iy}{(q,ab,ac,ad,bc,bd,cd;q)_\iy}\,,\quad
\frac{h_n}{h_0}=\frac{1-abcdq^{n-1}}{1-abcdq^{2n-1}}\,
\frac{(q,ab,ac,ad,bc,bd,cd;q)_n}{(abcd;q)_n}\,,
\label{eq:02.05}
\end{equation}
and the $x_k$ are the points $\thalf(eq^k+e^{-1}q^{-k})$ with $e$ any of the
$a,b,c,d$ of absolute value $>1$. The sum is over the
$k\in\Znonneg$ with $|eq^k|>1$; the sum does not occur if
moreover $|a|,|b|,|c|,|d|\le 1$.
The $\om_k$ for $e=a$ are the weights
\begin{equation}
\om_k=\frac{(a^{-2};q)_\iy}{(q,ab,ac,ad,a^{-1}b,a^{-1}c,a^{-1}d;q)_\iy}\,
\frac{(1-a^2 q^{2k})(a^2,ab,ac,ad;q)_k}
{(1-a^2)(q,qab^{-1},qac^{-1},qad^{-1};q)_k}\,
\left(\frac q{abcd}\right)^k.
\label{106}
\end{equation}

A more uniform way of writing the orthogonality relation \eqref{eq:02.03}
is by the contour integral
\begin{equation}
\frac1{2\pi i}
\oint_C p_n\bigl(\thalf(z+z^{-1})\bigr)\,p_m\bigl(\thalf(z+z^{-1})\bigr)\,
\frac{(z^2,z^{-2};q)_\iy}
{(az,az^{-1},bz,bz^{-1},cz,cz^{-1},dz,dz^{-1};q)_\iy}\,\frac{\dup z}z
=2h_n\de_{n,m},
\label{eq:02.06}
\end{equation}
where $C$ is the unit circle traversed in positive direction with suitable
deformations to separate the sequences of poles converging to zero
from the sequences of poles diverging to $\iy$.

The case $n=m=0$ of \eqref{eq:02.06} or \eqref{eq:02.03} is known as the
{\em Askey--Wilson integral}.
\paragraph{$q$-Difference equation}
\begin{equation}
\begin{split}
&LR_n=(q^{-n}-1)(1-q^{n-1}abcd)R_n,\\
&(Lf)(z):=A(z)f(qz)-\bigl(A(z)+A(z^{-1})\bigr)f(z)+A(z^{-1})f(q^{-1}z),\\
&A(z):=\frac{(1-az)(1-bz)(1-cz)(1-dz)}{(1-z^2)(1-qz^2)}\,.
\end{split}
\label{eq:02.07}
\end{equation}
\paragraph{Duality}
Assume that $\Re a>0$ and $q^{-1}abcd\in\CC\backslash (-\iy,0]$.
For $z\in\CC\backslash (-\iy,0]$ let $\sqrt{z}$ be such that
$\sqrt{z}>0$ if $z>0$.
Then the \emph{dual parameters} $\td a, \td b, \td c, \td d$ are well defined
by
\begin{equation}
\td a=\sqrt{q^{-1}abcd},
\quad\td a\td b=ab,
\quad\td a\td c=ac,
\quad\td a\td d=ad.
\label{102}
\end{equation}
Moreover, taking duals
of dual parameters once more, we recover $a,b,c,d$.
We can use the dual parameters to observe from \eqref{101} a 
duality of Askey--Wilson polynomials:
\begin{equation}
R_n(q^m a;a,b,c,d)
=R_m(q^n\td a;\td a,\td b,\td c,\td d)\quad
(m,n\in\ZZ_{\ge0}).
\label{103}
\end{equation}
\paragraph{Askey--Wilson functions}
By Watson's transformation \eqref{eq:01.48} Askey--Wilson polynomials
\eqref{101} can also be expressed as terminating ${}_8W_7$ functions.
These can be extended to \emph{Askey--Wilson functions}
\begin{multline}
\phi_\ga(z;a,b,c,d;q):=
\frac{(q\td a d^{-1}\ga z,q\td a d^{-1}\ga z^{-1};q)_\iy}
{(bc,q\td a ad^{-1}\ga,\td a^{-1}bc\ga,qd^{-1}z,qd^{-1}z^{-1};q)_\iy}
\\
\times{}_8W_7(\td a ad^{-1}\ga;az,az^{-1},
\td a\ga,\td b\ga,\td c\ga;q,\td a^{-1}bc\ga).
\label{120}
\end{multline}
Then, with $L$ given  by \eqref{eq:02.07},
\begin{equation}
L\phi_\ga=-(1-\td a\ga)(1-\td a\ga^{-1})\phi_\ga.
\end{equation}
For special $\ga$ we get Askey--Wilson polynomials:
\begin{equation}
\phi_{q^n\td a}(z;a,b,c,d;q)=\frac1{(bc,qad^{-1},qa^{-1}d^{-1};q)_\iy}\,
R_n(z;a,b,c,d\,|\,q)\qquad(n\in\ZZ_{\ge0}).
\end{equation}
See Koelink \& Stokman \cite{K-S2001}, who also define the
\emph{Askey--Wilson function transform} and give the inverse transform,
both transforms having \eqref{120} as a kernel.
\subsection{Continuous $q$-ultraspherical polynomials}
\paragraph{Definitions as finite Fourier series and as special
Askey--Wilson polynomial}
\begin{align}
C_n(\cos\tha;\be\,|\,q):=&
\sum_{k=0}^n\frac{(\be;q)_k(\be;q)_{n-k}}{(q;q)_k(q;q)_{n-k}}\,
\eup^{\iup(n-2k)\tha}
\label{eq:02.08}\\
=&\frac{(\be;q)_n}{(q;q)_n}\,
p_n(\cos\tha;\be^\half,q^\half\be^\half,-\be^\half,-q^\half\be^\half\,|\,q).
\label{eq:02.09}
\end{align}
\paragraph{Orthogonality relation}
($-1<\be<1$)
\begin{equation}
\frac1{2\pi}\int_0^\pi
C_n(\cos\tha;\be,q)\,C_m(\cos\tha;\be,q)\,
\left|\frac{(\eup^{2\iup\tha};q)_\iy}{(\be\eup^{2\iup\tha};q)_\iy}\right|^2
\dup\tha=\frac{(\be,q\be;q)_\iy}{(\be^2,q;q)_\iy}\,
\frac{1-\be}{1-\be q^n}\,\frac{(\be^2;q)_n}{(q;q)_n}\,\de_{n,m}.
\label{eq:02.10}
\end{equation}
\paragraph{$q$-Difference equation}
\begin{equation}
A(z)P_n(qz)-\bigl(A(z)+A(z^{-1})\bigr)P_n(z)+A(z^{-1})P_n(q^{-1}z)=
(q^{-n}-1)(1-q^n\be^2)P_n(z),
\label{eq:02.11}
\end{equation}
where $P_n(z)=C_n\bigl(\thalf(z+z^{-1});\be\,|\,q\bigr)$ and
$A(z)=(1-\be z^2)(1-q\be z^2)/\bigl((1-z^2)(1-qz^2)\bigr)$.
\paragraph{Generating function}
\begin{equation}
\frac{(\be \eup^{\iup\tha}z,\be \eup^{-\iup\tha}z;q)_\iy}
{(\eup^{\iup\tha}z,\eup^{-\iup\tha}z;q)_\iy}=
\sum_{n=0}^\iy C_n(\cos\tha;\be\,|\,q)\,z^n\qquad
(|z|<1,\;0\le\tha\le\pi,\;-1<\be<1).
\label{eq:02.12}
\end{equation}
\paragraph{Special case: the continuous $q$-Hermite polynomials}
\begin{equation}
H_n(x\,|\,q)=(q;q)_n\,C_n(x;0\,|\,q).
\label{eq:02.13}
\end{equation}
\paragraph{Special cases: the Chebyshev polynomials}
\begin{align}
C_n(\cos\tha;q\,|\,q)&=U_n(\cos\tha):=\frac{\sin((n+1)\tha)}{\sin\tha}\,,
\label{eq:02.14}\\
\lim_{\be\uparrow 1}\frac{(q;q)_n}{(\be;q)_n}\,
C_n(\cos\tha;\be\,|\,q)&=
T_n(\cos\tha):=\cos(n\tha)
\qquad(n>0).
\label{eq:02.15}
\end{align}
\subsection{$q$-Racah polynomials}
\paragraph{Definition as $q$-hypergeometric series}
($n=0,1,\ldots,N$)
\begin{equation}
R_n(q^{-y}+\ga\de q^{y+1};\al,\be,\ga,\de\,|\,q):=
\qhyp43{q^{-n},\al\be q^{n+1},q^{-y},\ga\de q^{y+1}}{q\al,q\be\de,q\ga}{q,q}
\qquad\mbox{($\al$, $\be\de$ or $\ga=q^{-N-1}$).}
\label{eq:02.16}
\end{equation}
\paragraph{Orthogonality relation}
\begin{equation}
\sum_{y=0}^N R_n(q^{-y}+\ga\de q^{y+1})\,R_m(q^{-y}+\ga\de q^{y+1})\,\om_y
=h_n\de_{n,m},
\label{eq:02.17}
\end{equation}
where $\om_y$ and $h_n$ can be explicitly given.
\subsection{Big $q$-Jacobi polynomials}
\paragraph{Definition as $q$-hypergeometric series}
\begin{equation}
P_n(x)=P_n(x;a,b,c;q):=\qhyp32{q^{-n},q^{n+1}ab,x}{qa,qc}{q,q}.
\label{eq:02.18}
\end{equation}
\paragraph{Orthogonality relation}
\begin{equation}
\int_{qc}^{qa}P_n(x)\,P_m(x)\,
\frac{(a^{-1}x,c^{-1}x;q)_\iy}{(x,bc^{-1}x;q)_\iy}\,\dup_qx=h_n\,\de_{n,m},
\quad(0<a<q^{-1},\;0<b<q^{-1},\;c<0),
\label{eq:02.19}
\end{equation}
where $h_n$ can be explicitly given.
\paragraph{$q$-Difference equation}
\begin{equation}
A(x)P_n(qx)-(A(x)+C(x))P_n(x)+C(x)P_n(q^{-1}x)=
(q^{-n}-1)(1-abq^{n+1})P_n(x),
\label{eq:02.20}
\end{equation}
where $A(x)=aq(x-1)(bx-c)/x^2$ and $C(x)=(x-qa)(x-qc)/x^2$.
\paragraph{Limit case: Jacobi polynomials $P_n^{(\al,\be)}(x)$}
\begin{equation}
\lim_{q\uparrow1} P_n(x;q^\al,q^\be,-q^{-1}d;q)
=\frac{n!}{(\al+1)_n}\,P_n^{(\al,\be)}\left(\frac{2x+d-1}{d+1}\right).
\label{eq:02.24}
\end{equation}
\paragraph{Special case: the little $q$-Jacobi polynomials}
\begin{align}
p_n(x;a,b;q)&=(-b)^{-n}q^{-\half n(n+1)}\,\frac{(qb;q)_n}{(qa;q)_n}\,
P_n(qbx;b,a,0;q)
\label{eq:02.21}\\
&={}_2\phi_1(q^{-n},q^{n+1}ab;qa;q,qx),
\label{eq:02.22}
\end{align}
which satisfy orthogonality relation (for $0<a<q^{-1}$ and $b<q^{-1}$)
\begin{equation}
\int_0^1 p_n(x;a,b;q)\,p_m(x;a,b;q)\,
\frac{(qx;q)_\iy}{(qbx;q)_\iy}\,x^{\log_q a}\,\dup_qx=
\frac{(q,qab;q)_\iy}{(qa,qb;q)_\iy}\,\frac{(1-q)(qa)^n}{1-abq^{2n+1}}
\frac{(q,qb;q)_n}{(qa,qab;q)_n}\,\de_{n,m}.
\label{eq:02.23}
\end{equation}
\paragraph{Limit case: Jackson's third $q$-Bessel function}
(see \eqref{eq:01.63} and \cite{20})
\begin{equation}
\lim_{N\to\iy}p_{N-n}(q^{N+k};q^\nu,b;q)=
\frac{(q;q)_\iy}{(q^{\nu+1};q)_\iy}\,q^{-\nu(n+k)}\,
J_\nu^{(3)}(2q^{\half(n+k)};q)\qquad(\nu>-1),
\label{eq:02.34}
\end{equation}
by which \eqref{eq:02.23} tends to the orthogonality relation for
$J_\nu^{(3)}(x;q)$:
\begin{equation}
\sum_{k=-\iy}^\iy J_\nu^{(3)}(2q^{\half(n+k)};q)\,
J_\nu^{(3)}(2q^{\half(m+k)};q)\,q^k=\de_{n,m} q^{-n}\qquad
(n,m\in\ZZ).
\label{eq:02.35}
\end{equation}
This orthogonlity relation can be equivalently written as a pair of a
$q$-Hankel transform and its inverse.
\subsection{$q$-Hahn polynomials}
\paragraph{Definition as $q$-hypergeometric series}
\begin{equation}
Q_n(x;\al,\be,N;q):=
\qhyp32{q^{-n},q^{n+1}\al\be,x}{q\al,q^{-N}}{q,q}
\qquad(n=0,1,\ldots,N).
\label{eq:02.25}
\end{equation}
\paragraph{Orthogonality relation}
\begin{equation}
\sum_{y=0}^N Q_n(q^{-y})\,Q_m(q^{-y})\,
\frac{(q\al,q^{-N};q)_y\,(q\al\be)^{-y}}{(q^{-N}\be^{-1},q;q)_y}
=h_n\de_{n,m},
\label{eq:02.26}
\end{equation}
where $h_n$ can be explicitly given.
\subsection{Stieltjes--Wigert polynomials}
\paragraph{Definition as $q$-hypergeometric series}
\begin{equation}
S_n(x;q)=\frac1{(q;q)_n}\,\qhyp11{q^{-n}}0{q,-q^{n+1}x}.
\label{eq:02.27}
\end{equation}
The orthogonality measure is not uniquely determined:
\begin{multline}
\int_0^\iy S_n(q^\half x;q)\,S_m(q^\half x;q)\,w(x)\,\dup x=
\frac1{q^n(q;q)_n}\,\de_{n,m},\quad
\mbox{where, for instance,}\\
w(x)=\frac{q^\half}{\log(q^{-1})\,(q,-q^\half x,-q^\half x^{-1};q)_\iy}
\quad{\rm or}\quad
\frac{q^\half}{\sqrt{2\pi\log(q^{-1})}}\,
\exp\left(-\frac{\log^2x}{2\log(q^{-1})}\right).
\label{eq:02.28}
\end{multline}
\subsection{Limits for $q\to-1$}
Bannai \& Ito \cite[pp.~271--273]{34} considered limits for $q\to-1$ of the
$q$-Racah polynomials \eqref{eq:02.16} (after suitable rescaling of the
$q$-Racah parameters).
The resulting polynomials, called
\emph{Bannai--Ito polynomials}, can be explicitly given as a sum of two
hypergeometric ${}_4F_3(1)$ polynomials, with the precise analytic form
depending on the parity of the degree $n$.
In \cite[(8.7)]{35} a version of the Bannai--Ito polynomials was obtained as a
limit for $q\to-1$ of the (rescaled) Askey--Wilson polynomials.
As a limit case of \eqref{eq:02.07} they are seen as eigenfunctions of
a Dunkl type (cf.~\S\ref{113}) difference-reflection operator
\cite[(9.3)--(9.7)]{35}.
Next, in \cite[(3.7)]{36} a version called \emph{continuous} Bannai--Ito polynomials
was shown to satisfy orthogonality relations on $(-\iy,\iy)$ with an explicit
weight function. In various papers by Vinet and coauthors limits for $q\to-1$ of
other families in the $q$-Askey scheme were considered. These can also be obtained
as limit cases of continuous or discrete Bannai--Ito polynomials.
In \cite{P-L-Z} an analogue of the Askey scheme is given for all such $q=-1$
orthogonal polynomials in the continuous case.
\subsection{Rahman--Wilson biorthogonal rational functions}
The following functions are rational in
their first argument:
\begin{equation}
R_n\left(\thalf(z+z^{-1});a,b,c,d,e\right):=
{}_{10}W_9(a/e;q/(be),q/(ce),q/(de),az,a/z,q^{n-1}abcd,q^{-n};q,q).
\label{eq:02.29}
\end{equation}
They satisfy the biorthogonality relation
\begin{equation}
\frac1{2\pi i}\oint_C 
R_n\left(\thalf(z+z^{-1});a,b,c,d,e\right)\,
R_m\left(\thalf(z+z^{-1});a,b,c,d,\frac q{abcde}\right)\,w(z)\,\frac{\dup z}z
=2h_n\de_{n,m},
\label{eq:02.30}
\end{equation}
where the contour $C$ is as in \eqref{eq:02.06}, and where
\begin{equation}
w(z)=
\frac{(z^2,z^{-2},abcdez,abcde/z;q)_\iy}
{(az,a/z,bz,b/z,cz,c/z,dz,d/z,ez,e/z;q)_\iy},
\label{eq:02.31}
\end{equation}
\begin{equation}
h_0=\frac{(bcde,acde,abde,abce,abcd;q)_\iy}
{(q,ab,ac,ad,ae,bc,bd,be,cd,ce,de;q)_\iy},
\label{eq:02.32}
\end{equation}
and $h_n/h_0$ can also be given explicitly.
For $ab=q^{-N}$, $n,m\in\{0,1,\ldots,N\}$ there is a related
discrete biorthogonality of the form
\begin{equation}
\sum_{k=0}^N
R_n\left(\thalf(aq^k+a^{-1}q^{-k});a,b,c,d,e\right)\,
R_m\left(\thalf(aq^k+a^{-1}q^{-k});a,b,c,d,\frac q{abcde}\right)\,
w_k=0\qquad(n\ne m).
\label{eq:02.33}
\end{equation}
See Rahman \cite{22} and Wilson \cite{21}.
\section{Orthogonal polynomials associated with root systems}
\label{sec:05}
\subsection{Macdonald polynomials for root system $A_{n-1}$}
\label{sec:05.03}
Reference for this subsection is Macdonald's book \cite[Ch.~VI]{09}].
Fix a positive integer $n$. We work with {\em partitions}
$\la=(\la_1,\ldots,\la_n)\in\ZZ^n$, where
$\la_1\ge\cdots\ge\la_n\ge0$.
Then $\la$ has \emph{weight} $|\la|:=\la_1+\cdots+\la_n$.
Note the special partition
\begin{equation}
\de:=(n-1,n-2,\ldots,0).
\label{119}
\end{equation}
On the set of such partitions consider the
{\em dominance partial ordering} $\le$ and the
{\em inclusion partial ordering} $\subseteq$:
\begin{align*}
\mu\le\la\quad&{\rm iff}\quad\mu_1+\cdots+\mu_i\le\la_1+\cdots+\la_i\quad
(i=1,\ldots,n);\\
\mu\subseteq\la\quad&{\rm iff}\quad\mu_i\le\la_i\quad(i=1,\ldots,n).
\end{align*}
The {\em monomials} are
$z^\al=z_1^{\al_1}\ldots z_n^{\al_n}$ ($\al_1,\ldots,\al_n\in\Znonneg$).
For $\la$ a partition and $S_n$ the symmetric group
the {\em symmetrized monomials} are defined by
\begin{equation}
m_\la(z):=\sum_{\al\in S_n\la} z^\al.
\label{eq:05.01}
\end{equation}
In integrals over the torus $T:=\{z\in \CC^n\,|\,|z_1|=\ldots=|z_n|=1\}$
write $\frac{\dup z}z:=\frac{\dup z_1}{z_1}\cdots\frac{\dup z_n}{z_n}$.
\paragraph{Definition}
For $\la$ a partition and $t\in(0,1)$
a {\em Macdonald polynomial}
is a polynomial of the form
\begin{equation}
P_\la(z)=P_\la(z;q,t)=m_\la(z)+
\sum_{\mu<\la;\,|\mu|=|\la|}u_{\la,\mu} m_\mu(z)
\label{115}
\end{equation}
such that for all $\mu<\la$ with $|\mu|=|\la|$
\begin{equation}
\frac1{(2\pi \iup)^n}\int_T P_\la(z)\,\overline{m_\mu(z)}\,\De(z)\,\frac{\dup z}z=0,
\label{116}
\end{equation}
where
\begin{equation}
\De(z):=|\De_+(z)|^2,\qquad
\De_+(z)=\De_+(z;q,t):=
\prod_{1\le i<j\le n}\frac{(z_iz_j^{-1};q)_\iy}
{(tz_iz_j^{-1};q)_\iy}\,.
\label{eq:05.03}
\end{equation}
\paragraph{Orthogonality relation}
For all partitions $\la,\mu$ we have
\begin{equation}
\frac1{(2\pi \iup)^n}\int_T P_\la(z)\,\overline{P_\mu(z)}\,\De(z)\,
\frac{\dup z}z=
n!\,\prod_{i<j}\frac{(q^{\la_i-\la_j}t^{j-i},q^{\la_i-\la_j+1}t^{j-i};q)_\iy}
{(q^{\la_i-\la_j}t^{j-i+1},q^{\la_i-\la_j+1}t^{j-i-1};q)_\iy}\,
\de_{\la,\mu}\,.
\label{eq:05.04}
\end{equation}
\paragraph{$q$-Difference equation}
\begin{equation}
\sum_{i=1}^n \prod_{j\ne i}\frac{tz_i-z_j}{z_i-z_j}\,\tau_{q,z_i}
P_\la(z;q,t)=\left(\sum_{i=1}^n q^{\la_i}t^{n-i}\right)P_\la(z;q,t),
\label{eq:05.05}
\end{equation}
where $\tau_{q,z_i}$ is the {\em $q$-shift operator}:
$\tau_{q,z_i}f(z_1,\ldots,z_n):=f(z_1,\ldots,qz_i,\ldots,z_n)$.
See \cite[Ch.~VI, \S3]{09} for the full system of
$q$-difference equations.
\paragraph{Special value}
\begin{equation}
P_\la(t^\de;q,t)=\prod_{i=1}^n t^{(i-1)\la_i}
\prod_{i<j}
\frac{(tq^{j-i};q)_{\la_i-\la_j}}{(q^{j-i};q)_{\la_i-\la_j}}\,.
\label{eq:05.06}
\end{equation}
\paragraph{Restriction of number of variables}
\begin{equation}
P_{\la_1,\la_2,\ldots,\la_{n-1},0}(z_1,\ldots,z_{n-1},0;q,t)=
P_{\la_1,\la_2,\ldots,\la_{n-1}}(z_1,\ldots,z_{n-1};q,t).
\label{eq:05.07}
\end{equation}
\paragraph{Homogeneity}
\begin{equation}
P_{\la_1,\ldots,\la_n}(z;q,t)=z_1\ldots z_n
P_{\la_1-1,\ldots,\la_n-1}(z;q,t)
\qquad(\la_n>0).
\label{eq:05.08}
\end{equation}
\paragraph{Self-duality}
Let $\la,\mu$ be partitions. Let
$q^\mu t^\de:=(q^{\mu_1}t^{n-1},q^{\mu_2}t^{n-2},\ldots,q^{\mu_n})$.
\begin{equation}
\frac{P_\la(q^\mu t^\de;q,t)}
{P_\la(t^\de;q,t)}
=
\frac{P_\mu(q^\la t^\de;q,t)}
{P_\mu(t^\de;q,t)}\,.
\label{eq:05.09}
\end{equation}
\paragraph{Interpolation Macdonald polynomials}(Sahi, Knop, Okounkov; see references
in \cite[\S5.1]{K2015})
The {\em interpolation Macdonald polynomial}
(or {\em shifted Macdonald polynomial})
$P_\la^{\rm ip}(z;q,t)$ is the unique $S_n$-invariant
polynomial of degree $|\la|$ with $z^\la$ having coefficient 1
such that
$P_\la^{\rm ip}(q^\mu t^\de;q,t)=0$
for each partition $\mu\ne\la$ with $|\mu|\le|\la|$.
They satisfy the \emph{extra vanishing property} that
$P_\la^{\rm ip}(q^\mu t^\de;q,t)=0$ if $\la\nsubseteq\mu$.
\paragraph{Binomial formula}(Okounkov; see reference in \cite[\S8.2]{K2015})
\begin{equation}
\frac{P_\la(z;q,t)}{P_\la(t^\de;q,t)}
=\sum_{\mu\subseteq\la}
\frac{P_\mu^{\rm ip}(q^\la t^\de;q,t)}{P_\mu^{\rm ip}(q^\mu t^\de;q,t}\,
\frac{P_\mu^{\rm ip}(z;q,t)}{P_\mu(t^\de;q,t)}\,.
\end{equation}
\paragraph{Special cases and limit relations}\quad\\
{\em Continuous $q$-ultraspherical polynomials} (see \eqref{eq:02.08}):
\begin{equation}
P_{m,n}(r \eup^{\iup\tha},r \eup^{-\iup\tha};q,t)=
\frac{(q;q)_{m-n}}{(t;q)_{m-n}}\,r^{m+n}\,C_{m-n}(\cos\tha;t\,|\,q).
\label{eq:05.10}
\end{equation}
{\em Symmetrized monomials} (see \eqref{eq:05.01}):
\begin{equation}
\lim_{t\uparrow1}P_\la(z;q,t)=m_\la(z).
\label{eq:05.11}
\end{equation}
{\em Schur functions}:
\begin{equation}
P_\la(z;q,q)=s_\la(z):=
\frac{\det(z_i^{\la_j+n-j})_{i,j=1,\ldots,n}}
{\det(z_i^{n-j})_{i,j=1,\ldots,n}}\,.
\label{eq:05.12}
\end{equation}
{\em Hall--Littlewood polynomials} (see \cite[Ch.~III]{09}):
\begin{equation}
P_\la(z;0,t)=P_\la(z;t).
\label{eq:05.14}
\end{equation}
{\em Jack polynomials} (see \cite[\S VI.10]{09}):
\begin{equation}
\lim_{q\uparrow1}P_\la(z;q,q^a)=P_\la^{(1/a)}(z).
\end{equation}
\paragraph{Algebraic definition of Macdonald polynomials}
Macdonald polynomials can also be defined
algebraically. We work now with partitions $\la=(\la_1,\la_2,\ldots)$
($\la_1\ge\la_2\ge\cdots\ge0$) such that $\la_j$ becomes eventually zero,
and with symmetric polynomials in arbitrarily many
variables $x_1,x_2,\ldots\,$, which can be canonically extended to
symmetric functions in infinitely many variables $x_1,x_2,\ldots\,$.
The \emph{length} $l(\la)$ of a partition $\la$ is the highest $j$ for which
$\la_j>0$.
The $r$th {\em power sum} $p_r$ and the symmetric functions $p_\la$
are formally defined by
\begin{equation}
p_r=\sum_{i\ge1}x_i^r,\qquad
p_\la=p_{\la_1}\,p_{\la_2}\ldots\,.
\label{eq:05.15}
\end{equation}
Put
\begin{equation}
z_\la:=\prod_{i\ge1}i^{m_i}\,m_i!\,,\quad
\mbox{where $m_i=m_i(\la)$ is the number of parts of $\la$ equal to $i$.}
\label{eq:05.16}
\end{equation}
Define an inner product $\lan\;,\;\ran_{q,t}$
on the space of symmetric functions such that
\begin{equation}
\lan p_\la,p_\mu\ran_{q,t}=\de_{\la,\mu}\,z_\la\,\prod_{i=1}^{l(\la)}
\frac{1-q^{\la_i}}{1-t^{\la_i}}\,.
\label{eq:05.17}
\end{equation}
For partitions $\la,\mu$ the partial ordering $\la\ge\mu$ means now that
$\sum_{j\ge1}\la_j=\sum_{j\ge1}\mu_j$ and
$\la_1+\cdots+\la_i\ge\mu_1+\cdots+\mu_i$ for all $i$.
The Macdonald polynomial $P_\la(x;q,t)$ can now be algebraically defined
as the unique symmetric function $P_\la$ of the form
$P_\la=m_\la+\sum_{\mu<\la}u_{\la,\mu}m_\mu$
such that
\begin{equation}
\lan P_\la,P_\mu\ran_{q.t}=0\qquad\mbox{if $\la\ne\mu$.}
\label{eq:05.18}
\end{equation}
If $l(\la)\le n$ then the newly defined $P_\la(x)$ with
$x_{n+1}=x_{n+2}=\ldots=0$ coincides with $P_\la(x;q,t)$ defined analytically,
and the new inner product is a constant multiple (depending on $n$) of the
old inner product.
\paragraph{Bilinear sum}
\begin{equation}
\sum_{\la}\frac1{\lan P_\la,P_\la\ran_{q,t}}\,
P_\la(x;q,t)P_\la(y;q,t)=\prod_{i,j\ge1}
\frac{(tx_iy_j;q)_\iy}{(x_iy_j;q)_\iy}.
\label{eq:05.19}
\end{equation}
\paragraph{Generalized Kostka numbers}
The {\em Kostka numbers} $K_{\la,\mu}$ occurring as expansion coefficients in
$s_\la=\sum_\mu K_{\la,\mu} m_\mu$ were generalized by Macdonald to
coefficients $K_{\la,\mu}(q,t)$ occurring in connection with Macdonald
polynomials, see \cite[\S VI.8]{09}.
Macdonald's conjecture that $K_{\la,\mu}(q,t)$ is a polynomial in $q$ and $t$
with coefficients in $\Znonneg$ was fully proved in \cite{10}.
\subsection{Orthogonal polynomials for general root systems}
\label{sec:05.04}
\paragraph{Constant term identitities}
Let $R$ be a reduced root system, $R^+$ the positive roots and $k\in\Zpos$.
Macdonald \cite{23} conjectured the second equality in 
\begin{equation}
\frac{\int_T \prod_{\al\in R^+}(\eup^{-\al};q)_k\,(q\eup^\al;q)_k\,\dup x}
{\int_T \dup x}=
\textup{CT}\Bigl(\prod_{\al\in R^+}
\prod_{i=1}^k(1-q^{i-1}\eup^{-\al})(1-q^i\eup^\al)
\Bigr)=\prod_{i=1}^n\qbinom {kd_i}kq,
\label{eq:05.20}
\end{equation}
where $T$ is a torus determined by $R$,
CT means the constant term in the Laurent expansion in $\eup^\al$,
and the $d_i$ are the degrees of the fundamental invariants of the
Weyl group of~$R$. For root system $A_{n-1}$ this becomes \eqref{eq:05.04}
for $\la=\mu=0$ and $t=q^k$.
The conjecture was extended for real $k>0$, for several
parameters $k$ (one for each root length), and for root system $BC_n$.
Next these conjectures became special cases of
Macdonald's conjectures \cite{15} for the
quadratic norms of Macdonald polynomials associated with root systems
(see next paragraph), and they were finally proved by Cherednik \cite{25}.

\paragraph{Macdonald polynomials for general root systems}
Macdonald \cite{15} introduced polynomials associated with an
arbitrary root system and depending on $q$, thus
extending the case of root system $A_{n-1}$ considered in
\S\ref{sec:05.03}.
These polynomials tend for $q\uparrow1$ to {\em Jacobi polynomials associated
with root systems}, which were introduced by Heckman and Opdam.
The paper \cite{15} ends with two conjectured explicit expressions:
(i) the quadratic norm of the polynomials (see \eqref{eq:05.04} for root system
$A_{n-1}$; the constant term conjecture is a special case of this conjecture);
(ii) evaluation of the polynomial for a special value
(see \eqref{eq:05.06} for root system $A_{n-1}$). In a lecture Macdonald also
formulated a duality conjecture, which becomes \eqref{eq:05.09} for $A_{n-1}$.
All these conjectures were proved by Cherednik \cite{25}, \cite{31}
by considering these Macdonald polynomials as
Weyl group symmetrizations of
non-invariant polynomials which are related to
a DAHA (\emph{double affine Hecke algebra}).
\subsection{Koornwinder polynomials}
Macdonald polynomials for root system $BC_n$ form a three-parameter
family which can be extended to the five-parameter
\emph{Koornwinder polynomials} (or \emph{Macdonald--Koornwinder polynomials})
\cite{16}. For the definition assume that $a,b,c,d$ are four reals, or two reals
and one pair of complex
conjugates, or two pairs of complex conjugates. Also assume that
pairwise products of $a,b,c,d$ are not equal to~1 and that
$|a|,|b|,|c|,|d|\le1$. Let $0<t<1$.
Let the group $W_n$ be the semidirect product of the symmetric group $S_n$
and the group $(\ZZ_2)^n$. This group acts on the \emph{Laurent monomials}
$z^\al=z_1^{\al_1}\ldots z_n^{\al_n}$ ($\al_1,\ldots,\al_n\in\ZZ$) by permutations
and sign changes of $\al_1,\ldots,\al_n$. We will work again with partitions
$\la=(\la_1,\ldots,\la_n)$, see \S\ref{sec:05.03}.
For $\la$ a partition the
\emph{symmetrized Laurent monomials} are defined by
\begin{equation}
\wt m_\la(z):=\sum_{\al\in W_n\la} z^\al.
\end{equation}
\paragraph{Definition}
For $\la$ a partition
a \emph{Koornwinder polynomial}
is a Laurent polynomial of the form
\begin{equation}
P_\la(z)=P_\la(z;q,t;a,b,c,d)=\wt m_\la(z)+
\sum_{\mu<\la}u_{\la,\mu} \wt m_\mu(z)
\end{equation}
such that for all $\mu<\la$
\begin{equation}
\frac1{(2\pi\iup)^n}\int_T P_\la(z)\,\wt m_\mu(z)\,\De(z)\,\frac{\dup z}z=0,
\end{equation}
where $\De(z)=|\De_+(z)|^2$ and
\begin{equation}
\De_+(z)=\De_+(z;q,t;a,b,c,d):=
\prod_{j=1}^n\frac{(z_j^2;q)_\iy}
{(az_j,bz_j,cz_j,dz_j;q)_\iy}\,
\prod_{1\le i<j\le n}\frac{(z_iz_j,z_iz_j^{-1};q)_\iy}
{(tz_iz_j,tz_iz_j^{-1};q)_\iy}.
\end{equation}
Then $P_\la(z)$ is real-valued and
\begin{equation}
\frac1{(2\pi\iup)^n}\int_T P_\la(z)\,P_\mu(z)\,\De(z)\,
\frac{\dup z}z=h_\la \de_{\la,\mu}
\label{117}
\end{equation}
for all partitions $\la,\mu$, where $h_\la$ is explicitly given in
\cite[(2.21)--(2.26)]{37}. In particular, foe $\la=\mu=0$ we have Gustafson's
\cite{24} formula
\begin{equation}
h_0=2^n n!\prod_{j=1}^n\frac{(t,t^{n+j-2}abcd;q)_\iy}
{(t^j,q,abt^{j-1},act^{j-1},\ldots,cdt^{j-1};q)_\iy}\,.
\label{eq:05.21}
\end{equation}

For $n=1$ the Koornwinder polynomials become the Askey--Wilson polynomials
\eqref{101}, multiplied
by a constant factor such that the coefficient of $z^n$ is 1.
\paragraph{$q$-Difference equation}
There is a second order $q$-difference operator $D$, explicitly given in
\cite[(5.4)]{16}, such that
\begin{equation}
DP_\la=c_\la P_\la,\qquad
c_\la=\sum_{j=1}^n\big(q^{-1}abcd\,t^{2n-j-1}(q^{\la_j}-1)
+t^{j-1}(q^{-\la_j}-1)\big).
\end{equation}
\paragraph{Duality}
Let $\de$ be the special partition \eqref{119} and let
$\td a,\td b,\td c,\td d$ be the dual parameters \eqref{102}.
For partitions $\la,\nu$ there is a duality formula including for $n=1$ the
Askey--Wilson case \eqref{103}:
\begin{equation}
\frac{P_\la(q^\nu t^\de a;q,t;a,b,c,d)}
{P_\la(t^\de a;q,t;a,b,c,d)}
=\frac{P_\nu(q^\la t^\de \td a;q,t;\td a,\td b,\td c,\td d)}
{P_\nu(t^\de \td a;q,t;\td a,\td b,\td c,\td d)}\,.
\label{118}
\end{equation}
The denominator on the left-hand side can be explicitly evaluated,
see \cite[(4.15)]{K2015}.

The duality \eqref{118}, the denominator evaluation in \eqref{118}, and the
evaluation of $h_\la$ in \eqref{117}
were all conjectured by Macdonald
(1991, unpublished). The duality conjecture was proved by Sahi \cite{32}.
By results of van Diejen \cite{33} this implied the other two conjectures.
\paragraph{$BC_n$-type interpolation Macdonald polynomials}
(Okounkov, see reference in \cite[\S5.3]{K2015})
The \emph{$BC_n$-type interpolation Macdonald polynomial}
$P_\la^{\rm ip}(z;q,t,a)$ is the unique $W_n$-invariant Laurent
polynomial of degree $|\la|$ with $z^\la$ having coefficient 1
such that $P_\la^{\rm ip}(q^\mu t^\de a;q,t,a)=0$
for each partition $\mu\ne\la$ with $|\mu|\le|\la|$.
They satisfy the extra vanishing property that
$P_\la^{\rm ip}(q^\mu t^\de a;q,t,a)=0$ if $\la\nsubseteq\mu$.
\paragraph{Binomial formula}(Okounkov; see reference in \cite[\S8.1]{K2015})
\begin{equation}
\frac{P_\la(z;q,t;a,b,c,d)}{P_\la(t^\de a;q,t;a,b,c,d)}
=\sum_{\mu\subseteq\la}
\frac{P_\mu^{\rm ip}(q^\la t^\de \td a;q,t,\td a)}
{P_\mu^{\rm ip}(q^\mu t^\de \td a;q,t,\td a)}\,
\frac{P_\mu^{\rm ip}(z;q,t,a)}{P_\mu(t^\de a;q,t;a,b,c,d)}\,.
\end{equation}
\subsection{Non-symmetric orthogonal polynomials}
\label{113}
Dunkl \cite{30} introduced \emph{Dunkl operators}:
certain differential-reflection operators
associated with a Weyl group (or a Coxeter group).
Cherednik \cite{25} introduced much more involved analogues of the Dunkl operators
in the context of the polynomial representation of a DAHA associated with a
reduced root system. His operators are certain $q$-difference-reflection operators.
For the non-reduced root system $BC_n$ these operators were introduced by Sahi
\cite{32}.
In all Dunkl type situations the operators commute and their joint eigenfunctions
are so-called \emph{non-symmetric} special functions, which yield more familiar
symmetric special functions by Weyl group symmetrization.
For root system $A_1$ or $BC_1$ the non-symmetric special functions can usually be
written as a sum of two terms, the first term being the symmetric special function
and the second term the anti-symmetric function, which can still be expressed in
terms of the symmetric special function, but with shifted parameters.

The non-symmetric Askey--Wilson polynomials
$E_{-n}(z;a,b,c,d\,|\,q)$ ($n\ge1$) and $E_n(z;a,b,c,d\,|\,q)$ ($n\ge0$),
given in \cite{N-S}, \cite{K2007}, are multiples of
\begin{equation}
R_n(z;a,b,c,d\,|\,q)
-Cz^{-1}(1-az)(1-bz)R_{n-1}(z;qa,qb,c,d\,|\,q)
\label{104}
\end{equation}
with
\begin{equation}
C=\frac{q^{1-n}(1-q^{n-1}abcd)(1-q^n ab)}{b(1-ab)(1-qab)(1-ac)(1-ad)}\quad\mbox{and}
\quad
\frac{aq^{1-n}(1-q^n)(1-q^{n-1}cd)}{(1-ab)(1-qab)(1-ac)(1-ad)}\,,
\label{105}
\end{equation}
respectively. They are Laurent polynomials, $E_{-n}(z)$ being spanned by
$z^{-n},\ldots,z^{n-1}$ and $E_n(z)$ by $z^{-n},\ldots,z^n$.

Define an inner product $\lan\,\cdot\,,\,\cdot\,\ran_{a,b,c,d;q}$
on the space of polynomials such that the left-hand side of \eqref{eq:02.03}
equals $h_0\lan p_n,p_m\ran_{a,b,c,d;q}$.
On the space of Laurent polynomials of the form
\begin{equation}
f(z)=f_1\big(\thalf(z+z^{-1})\big)+z^{-1}(1-az)(1-bz)f_2\big(\thalf(z+z^{-1})\big)
\quad\mbox{($f_1,f_2$ ordinary polynomials)},
\label{107}
\end{equation}
define an inner product
\begin{equation}
\lan g,h\ran:=\lan g_1,h_1\ran_{a,b,c,d;q}
-ab\,\frac{(1-ab)(1-qab)(1-ac)(1-ad)(1-bc)(1-bd)}{(1-abcd)(1-qabcd)}\,
\lan g_2,h_2\ran_{qa,qb,c,d;q}\,,
\label{108}
\end{equation}
In fact, we have rewritten Laurent polynomials as vector-valued ordinary polynomials,
and we have defined an inner product on these vector-valued polynomials.
The Laurent polynomials $E_n(z)$ ($n\in\ZZ$) are orthogonal with respect
to the inner product \eqref{108}, see \cite{K-B2011}. If, in addition to the
conditions
for \eqref{eq:02.03}, $a$ and $b$ are real with opposite signs then the inner
product \eqref{108} is positive definite. 
\section{Elliptic hypergeometric series}
\label{sec:04}
See Gasper \& Rahman \cite[Ch.~11]{01} and references given there.\\
Let $p,q\in\CC$, $|p|,|q|<1$. Define a {\em modified Jacobi theta function}
by
\begin{equation}
\tha(x;p):=(x,p/x;p)_\iy\qquad(x\ne0),
\label{eq:04.01}
\end{equation}
and the {\em elliptic shifted factorial} by
\begin{equation}
(a;q,p)_k:=\tha(a;p)\tha(aq;p)\ldots\tha(aq^{k-1};p)\quad(k\in\Zpos);
\qquad
(a;q,p)_0:=1,
\label{eq:04.02}
\end{equation}
\begin{equation}
(a_1,\ldots,a_r;q,p)_k:=(a_1;q,p)_k\ldots(a_r;q,p)_k,
\label{eq:04.03}
\end{equation}
where $a,a_1,\ldots,a_r\ne0$. For $q=\eup^{2\pi\iup\si}$,
$p=\eup^{2\pi\iup\tau}$
($\Im\tau>0$) and $a\in\CC$ we have
\begin{equation}
\frac{\tha(a \eup^{2\pi\iup\si(x+\si^{-1})};\eup^{2\pi\iup\tau})}
{\tha(a \eup^{2\pi\iup\si x};\eup^{2\pi\iup\tau})}=1,\qquad
\frac{\tha(a \eup^{2\pi\iup\si(x+\tau\si^{-1})};\eup^{2\pi\iup\tau})}
{\tha(a \eup^{2\pi\iup\si x};\eup^{2\pi\iup\tau})}=-a^{-1}q^{-x}.
\label{eq:04.10}
\end{equation}

A series $\sum_{k=0}^\iy c_k$ with $c_{k+1}/c_k$ being an elliptic
(i.e.\ doubly periodic meromorphic) function of $k$ considered as a complex
variable, is called an {\em elliptic hypergeometric series}.
In particular, define
the ${}_{r}E_{r-1}$ {\em theta hypergeometric series} as the formal series
\begin{equation}
{}_{r}E_{r-1}(a_1,\ldots,a_r;b_1,\ldots,b_{r-1};q,p;z):=
\sum_{k=0}^\iy\frac{(a_1,\ldots,a_r;q,p)_k}{(b_1,\ldots,b_{r-1};q,p)_k}\,
\frac{z^k}{(q;q,p)_k}\,.
\label{eq:04.04}
\end{equation}
It has $g(k):=c_{k+1}/c_k$ with
\[
g(x)=\frac{z\,\tha(a_1q^x;p)\ldots\tha(a_r q^x;p)}
{\tha(q^{x+1};p)\,\tha(b_1q^x;p)\ldots\tha(b_{r-1}q^x;p)}\,.
\]
By \eqref{eq:04.10} 
$g(x)$ is an elliptic function with periods
$\si^{-1}$ and $\tau\si^{-1}$ ($q=e^{2\pi i\si}$, $p=e^{2\pi i\tau}$) if
the balancing condition
$a_1\ldots a_r=q b_1\ldots b_{r-1}$ is satisfied.

The ${}_rV_{r-1}$ {\em very-well-poised theta hypergeometric series}
(a special ${}_rE_{r-1}$) is defined, in case of argument 1, as:
\begin{equation}
{}_rV_{r-1}(a_1;a_6,\ldots,a_r;q,p):=
\sum_{k=0}^\iy\frac{\tha(a_1q^{2k};p)}{\tha(a_1;p)}\,
\frac{(a_1,a_6,\ldots,a_r;q,p)_k}
{(qa_1/a_6,\ldots,qa_1/a_r;q,p)_k}\,\frac{q^k}{(q;q,p)_k}\,.
\label{eq:04.05}
\end{equation}
The series is called {\em balanced} if
$a_6^2\ldots a_r^2=a_1^{r-6}q^{r-4}$. The series terminates if, for instance,
$a_r=q^{-n}$.
\paragraph{Elliptic analogue of Jackson's ${}_8W_7$ summation formula
\eqref{eq:01.47}}
\begin{equation}
{}_{10}V_9(a;b,c,d,q^{n+1}a^2/(bcd),q^{-n};q,p)
=\frac{(qa,qa/(bc),qa/(bd),qa/(cd);q,p)_n}{(qa/b,qa/c,qa/d,qa/(bcd);q,p)_n}\,.
\label{eq:04.06}
\end{equation}
\paragraph{Elliptic analogue of Bailey's ${}_{10}W_9$ transformation formula
\eqref{eq:01.58}}
\begin{multline}
{}_{12}V_{11}\left(a;b,c,d,e,f,\frac{q^{n+2}a^3}{bcdef},q^{-n};q,p\right)
=\frac{(qa,qa/(ef),(qa)^2/(bcde),(qa)^2/(bcdf);q,p)_n}
{(qa/e,qa/f,(qa)^2/(bcdef),(qa)^2/(bcd);q,p)_n}\\
\times
{}_{12}V_{11}\left(\frac{qa^2}{bcd};\frac{qa}{cd},\frac{qa}{bd},\frac{qa}{bc},
e,f,\frac{q^{n+2}a^3}{bcdef},q^{-n};q,p\right).
\label{eq:04.07}
\end{multline}
Suitable ${}_{12}V_{11}$ functions satisfy a discrete biorthogonality
relation which is an elliptic analogue of
\eqref{eq:02.33}.
\paragraph{Ruijsenaars' elliptic gamma function}
(see Ruijsenaars \cite{26})
\begin{equation}
\Ga(z;q,p):=\prod_{j,k=0}^\iy\frac{1-z^{-1}q^{j+1}p^{k+1}}{1-zq^jp^k}\,,
\label{eq:04.08}
\end{equation}
which is symmetric in $p$ and $q$.
Then
\begin{equation}
\Ga(qz;q,p)=\tha(z;p)\,\Ga(z;q,p),\qquad
\Ga(q^nz;q,p)=(z;q,p)_n\,\Ga(z;q,p).
\label{eq:04.09}
\end{equation}
\section{Applications}
\label{sec:03}
\subsection{Quantum groups}
A specific quantum group is usually a Hopf algebra which is a
$q$-deformation of the Hopf algebra of functions on a specific Lie group
or, dually, of a universal enveloping algebra (viewed as Hopf algebra) of a
Lie algebra. The general philosophy is that representations of
the Lie group or Lie algebra also deform to representations of the quantum
group, and that special functions associated with the representations
in the classical case deform to $q$-special functions associated with
the representations in the quantum case. Sometimes this is straightforward,
but often new subtle phenomena occur.

The representation theoretic objects
which may be explicitly written in terms of $q$-special functions
include matrix elements of representations with respect to specific
bases (in particular spherical elements), Clebsch-Gordan coefficients
and Racah coefficients. Many one-variable $q$-hypergeometric functions
have found interpretation in some way in connection with a quantum
analogue of a three-dimensional Lie group (generically the Lie group
$SL(2,\CC)$ and its real forms).
Classical by now are: little $q$-Jacobi polynomials interpreted as
matrix elements of irreducible representations of $SU_q(2)$ with respect
to the standard basis; Askey--Wilson polynomials similarly interpreted
with respect to a certain basis not coming from a quantum subgroup;
Jackson's third $q$-Bessel functions as matrix elements of irreducible
representations of $E_q(2)$; $q$-Hahn polynomials and $q$-Racah polynomials
interpreted as Clebsch-Gordan coefficients and Racah coefficients,
respectively, for $SU_q(2)$.
See for instance Vilenkin \& Klimyk \cite{07} and Koelink \cite{27}

Further developments include: Macdonald polynomials as spherical elements
on quantum analogues of compact Riemannian symmetric spaces;
$q$-analogues of Jacobi functions as matrix elements of irreducble unitary
representations of $SU_q(1,1)$; Askey--Wilson polynomials as
matrix elements of representations of the $SU(2)$ dynamical quantum group;
an interpretation of discrete ${}_{12}V_{11}$ biorthogonality relations
on the elliptic $U(2)$ quantum group (see \cite{28}).

Since the $q$-deformed Hopf algebras are usually presented by generators
and relations, identities for $q$-special functions involving
non-commuting variables satisfying simple relations are important
for further interpretations of $q$-special functions in quantum groups,
for instance:
\paragraph{$q$-Binomial formula with $q$-commuting variables}
\begin{equation}
(x+y)^n=\sum_{k=0}^n\qbinom nkq y^{n-k} x^k\qquad(xy=qyx).
\label{eq:03.01}
\end{equation}
\paragraph{Functional equations for $q$-exponentials with $xy=qyx$}
\begin{align}
e_q(x+y)=e_q(y)e_q(x),&\qquad
E_q(x+y)=E_q(x)E_q(y),
\label{eq:03.02}\\
e_q(x+y-yx)=e_q(x)e_q(y),&\qquad
E_q(x+y+yx)=E_q(y)E_q(x).
\label{eq:03.03}
\end{align}
See \cite{06} for further reading.
\subsection{Various algebraic settings}
\paragraph{Classical groups over finite fields (Chevalley groups)}\quad\\
$q$-Hahn polynomials and various kinds of $q$-Krawtchouk polynomials have
interpretations as spherical and intertwining functions on
classical groups ($GL_n$, $SO_n$, $Sp_n$) over a finite field $F_q$
with respect to suitable subgroups, see
Stanton \cite{11}.
\paragraph{Affine Kac-Moody algebras}
(see \cite{13})\\
The Rogers-Ramanujan identities \eqref{eq:01.49}, \eqref{eq:01.50}
and some of their generalizations were
interpreted in the context of characters of representations of
the simplest affine Kac-Moody algebra $A_1^{(1)}$.
Macdonald's generalization \cite{12} of Weyl's denominator formula to affine root
systems has an interpretation as an identity for the denominator
of the character of a representation of an affine Kac-Moody algebra.
\subsection{Partitions of positive integers}
\label{122}
Let $n$ be a positive integer, $p(n)$ the number of partitions of $n$,
$p_N(n)$ the number of partitions of $n$ into parts $\le N$,
$p_{\rm dist}(n)$ the number of partitions of $n$ into distinct parts,
and $p_{\rm odd}(n)$ the number of partitions of $n$ into odd parts.
Then Euler observed:
\begin{align}
\frac1{(q;q)_\iy}=\sum_{n=0}^\iy p(n) q^n,\qquad
&\frac1{(q;q)_N}=\sum_{n=0}^\iy p_N(n) q^n,
\label{eq:03.05}\\
(-q;q)_\iy=\sum_{n=0}^\iy p_{\rm dist}(n) q^n,\qquad
&\frac1{(q;q^2)_\iy}=\sum_{n=0}^\iy p_{\rm odd}(n) q^n,
\label{eq:03.07}
\end{align}
and
\begin{equation}
(-q;q)_\iy=\frac1{(q;q^2)_\iy}\,,\qquad
p_{\rm dist}(n)=p_{\rm odd}(n).
\label{eq:03.08}
\end{equation}

The Rogers--Ramanujan identity \eqref{eq:01.49} has the following
partition theoretic interpretation:
The number of partitions of $n$ with parts differing at least 2 equals
the number of partitions of $n$ into parts congruent to 1 or 4$\pmod 5$.
Similarly, \eqref{eq:01.50} yields:
The number of partitions of $n$ with parts larger than 1 and
differing at least 2 equals
the number of partitions of $n$ into parts congruent to 2 or 3$\pmod 5$.

The left-hand sides of the Rogers--Ramanujan identities \eqref{eq:01.49} and
\eqref{eq:01.50} have interpretations in the {\em Hard Hexagon Model},
see \cite{18}. Much further work has been done on
Rogers--Ramanujan type identities in connection with more general models
in statistical mechanics. So-called {\em fermionic expressions} do occur.
%


\begin{thebibliography}{99}
%
\bibitem{19}
G. E. Andrews,
{\em $q$-Series: their development and application in analysis, number theory, combinatorics, physics, and computer algebra},
CBMS Regional Conference Series in Mathematics 66,
Amer. Math. Soc., 1986.
%
\bibitem{14}
G. E. Andrews, R. Askey and R. Roy,
{\em Special functions},
Cambridge University Press, 1999.
%
\bibitem{08}
G. E. Andrews and K. Eriksson,
{\em Integer partitions}.
Cambridge University Press, 2004.
%
\bibitem{29}
R. Askey and J. Wilson,
{\em Some basic hypergeometric orthogonal polynomials that generalize Jacobi 
polynomials},
Mem. Amer. Math. Soc. (1985), no. 319.
%
\bibitem{34}
E. Bannai and T. Ito,
\emph{Algebraic combinatorics. I: Association schemes},
Benjamin-Cummings, 1984.
%
\bibitem{18}
R. J. Baxter,
{\em Exactly solved models in statistical mechanics},
Academic Press, 1982.
%
\bibitem{25}
I. Cherednik,
{\em Double affine Hecke algebras and Macdonald's conjectures},
Ann. of Math. (2) 141 (1995), 191--216.
%
\bibitem{31}
I. Cherednik
{\em Macdonald's evaluation conjectures and difference Fourier transform},
Invent. Math. 122 (1995), 119--145; erratum: Invent. Math. 125 (1996), p. 391.
%
\bibitem{33}
J. F. van Diejen,
{\em Self-dual Koornwinder--Macdonald polynomials},
Invent. Math. 126 (1996), 319--339.
%
\bibitem{30}
C. F. Dunkl,
{\em Differential-difference operators associated to reflection groups},
Trans. Amer. Math. Soc. 311 (1989), 167--183.
%
\bibitem{01}
G. Gasper and M. Rahman,
{\em Basic hypergeometric series}, 2nd edn.
Cambridge University Press, 2004.
%
\bibitem{36}
V. X. Genest, L. Vinet and A. Zhedanov,
\emph{The non-symmetric Wilson polynomials are the Bannai--Ito polynomials},
Proc. Amer. Math. Soc. 144 (2016), 5217--5226.
%
\bibitem{G-L-Z}
Ya. I. Granovskii, I. M. Lutzenko and A. S. Zhedanov, 
\emph{Mutual integrability, quadratic algebras and dynamical symmetry},
Ann. Physics 217 (1992), 1--20.
%
\bibitem{24}
R. A. Gustafson,
{\em A generalization of Selberg's beta integral},
Bull. Amer. Math. Soc. (N.S.) 22 (1990), 97--105.
%
\bibitem{10}
M. Haiman,
{\em Hilbert schemes, polygraphs and the Macdonald positivity conjecture},
J. Amer. Math. Soc. 14 (2001), 941--1006.
%
\bibitem{02}
R. Koekoek, P.~A. Lesky and R.~F. Swarttouw,
{\em Hypergeometric orthogonal polynomials and their $q$-analogues},
Springer-Verlag, Berlin, 2010.
%
\bibitem{27}
H. T. Koelink,
{\em Askey--Wilson polynomials and the quantum $SU(2)$ group: survey and 
applications},
Acta Appl. Math. 44 (1996), 295--352.
%
\bibitem{K-S2001}
E. Koelink and J. V. Stokman,
\emph{The Askey--Wilson function transform},
Internat. Math. Res. Notices (2001), No. 22, 1203--1227.
%
\bibitem{28}
E. Koelink, Y. van Norden and H. Rosengren,
{\em Elliptic $U(2)$ quantum group and elliptic hypergeometric series},
Comm. Math. Phys. 245 (2004), 519--537.
%
\bibitem{16}
T. H. Koornwinder,
{\em Askey--Wilson polynomials for root systems of type $BC$},
in: 
{\em Hypergeometric functions on domains of positivity, Jack polynomials,
and applications}, D. St. P. Richards (ed.),
Contemp. Math. 138, Amer. Math. Soc., 1992, pp. 189--204.
%
\bibitem{04}
T. H. Koornwinder,
{\em Compact quantum groups and $q$-special functions},
in: {\em Representations of Lie groups and quantum groups}, 
V. Baldoni and M. A. Picardello (eds.), 
Pitman Research Notes in Mathematics Series 311,
Longman Scientific \& Technical, 1994, pp 46--128;
Sections 3 and 4 revised in
{\tt arXiv:math.CA/9403216v2}, 2013.
%
\bibitem{06}
T. H. Koornwinder,
{\em Special functions and $q$-commuting variables},
in: 
{\em Special functions, $q$-series and related topics},
M. E. H. Ismail, D. R. Masson and M. Rahman (eds.),
Fields Institute Communications 14, Amer. Math. Soc., 1997,
pp. 131--166;
{\tt arXiv:q-alg/9608008}.
%
\bibitem{K2007}
T. H. Koornwinder, The relationship between Zhedanov's algebra AW(3) and the double
affine Hecke algebra in the rank one case,
SIGMA 3 (2007), 063, 15 pp.; \texttt{arXiv:math/0612730v4}.
%
\bibitem{K2015}
T. H. Koornwinder,
\emph{Okounkov's $BC$-type interpolation Macdonald polynomials and their $q=1$ 
limit}, S\'em. Lothar. Combin. B72a (2015), 27 pp.; revised in
\texttt{arXiv:1408.5993v5}, 2015. 
%
\bibitem{K2022a}
T.~H. Koornwinder,
\emph{Charting the $q$-Askey scheme},
in: \emph{Hypergeometry, Integrability and Lie Theory},
Contemporary Mathematics 780 (2022), 79--94;  
\texttt{arXiv:2108.03858}
%
\bibitem{K2022b}
T. H. Koornwinder,
\emph{Charting the q-Askey scheme. II. The q-Zhedanov scheme},
\texttt{arXiv:2209.07995}, 2022. 
%
\bibitem{K-B2011}
T. H. Koornwinder and F. Bouzeffour,
\emph{Nonsymmetric Askey-Wilson polynomials as vector-valued polynomials},
Appl. Anal. 90 (2011), 731--746; revised in \texttt{arXiv:1006.1140v3}, 2018.
%
\bibitem{20}
T. H. Koornwinder and R. F. Swarttouw,
{\sl On $q$-Analogues of the Fourier and Hankel transforms},
Trans. Amer. Math. Soc.
333
(1992),
445--461.
%
\bibitem{13}
J. Lepowsky,
Affine Lie algebras and combinatorial identities,
in: {\em Lie algebras and related topics}, D. J. Winter (ed.),
Lecture Notes in Math. 933, Springer, 1982, pp. 130--156.
%
\bibitem{12}
I. G. Macdonald,
{\em Affine root systems and Dedekind's $\eta$-function},
Invent. Math. 15 (1972), 91--143.
%
\bibitem{23}
I. G. Macdonald, 
{\em Some conjectures for root systems},
SIAM J. Math. Anal. 13 (1982), 988--1007.
%
\bibitem{09}
I. G. Macdonald,
{\em Symmetric functions and Hall polynomials}, 2nd edn.
Clarendon Press, Oxford, 1995.
%
\bibitem{15}
I. G. Macdonald, 
{\em Orthogonal polynomials associated with root systems},
S\'em. Lothar. Combin. 45 (2000/01), Art. B45a;
{\tt arXiv:math.QA/0011046} (publication of a manuscript from 1987).
%
\bibitem{17}
I. G. Macdonald,
{\em Affine Hecke algebras and orthogonal polynomials},
Cambridge University Press, 2003.
%
\bibitem{N-S}
M. Noumi and J. V. Stokman,
\emph{Askey--Wilson polynomials: an affine Hecke algebraic approach}, in:
\emph{Laredo Lectures on Orthogonal Polynomials and Special Functions},
Nova Sci. Publ., Hauppauge, NY (2004), pp. 111--144;
\texttt{arXiv:math/0001033}.
%
\bibitem{P-L-Z}
J. Pelletier, L. Vinet and A. Zhedanov,
\emph{Continuous $-1$ hypergeometric orthogonal polynomials},
\texttt{arXiv:2209.10727}, 2022.
%
\bibitem{22}
M. Rahman,
{\em An integral representation of a ${}_{10}\phi_9$ and continuous biorthogonal 
${}_{10}\phi_9$ rational functions},
Canad. J. Math. 38 (1986), 605--618.
%
\bibitem{26}
S. N. M. Ruijsenaars,
{\em Special functions defined by analytic difference equations},
in: {\em Special functions 2000: current perspective and future directions},
J. Bustoz, M.~E.~H. Ismail and S.~K. Suslov (eds.),
NATO Sci. Ser. II Math. Phys. Chem. 30,
Kluwer, 2001, pp. 281--333.
%
\bibitem{32}
S. Sahi,
{\em Nonsymmetric Koornwinder polynomials and duality},
Ann. of Math. (2) 150 (1999), 267--282.
%
\bibitem{11}
D. Stanton,
{\em Orthogonal polynomials and Chevalley groups},
in: {\em Special functions: group theoretical aspects and applications},
R. A. Askey, T. H. Koornwinder and W. Schempp (eds.),
Reidel, Dordrecht, 1984, pp 87--128.
%
\bibitem{37}
J. V. Stokman,
\emph{On $BC$ type basic hypergeometric orthogonal polynomials},
Trans. Amer. Math. Soc. 352 (2000), 1527--1579.
%
\bibitem{05}
S. K. Suslov
{\em An introduction to basic Fourier series},
Kluwer, Dordrecht, 2003.
%
\bibitem{35}
S. Tsujimoto, L. Vinet and A. Zhedanov,
\emph{Dunkl shift operators and Bannai-Ito polynomials},
Adv. Math. 229 (2012), 2123--2158
%
\bibitem{03}
J. Van der Jeugt and K. Srinivasa Rao,
{\em Invariance groups of transformations of basic hypergeometric series},
J. Math. Phys. 40 (1999), 6692--6700.
%
\bibitem{VS2021}
L. Verde-Star,
\emph{A unified construction of all the hypergeometric and basic hypergeometric
families of orthogonal polynomial sequences},
Linear Algebra Appl. 627 (2021), 242--274.
%
\bibitem{V-Z2016}
L. Vinet and  A. Zhedanov,
\emph{Hypergeometric orthogonal polynomials with respect to Newtonian basis},
SIGMA 12 (2016), paper 048, 14 pp.;
\texttt{arXiv:1602.02724}.
%
\bibitem{07}
N. J. Vilenkin and A. U. Klimyk,
{\em Representation of Lie groups and special functions, Vol. 3},
Kluwer, Dordrecht, 1992.
%
\bibitem{21}
J. A. Wilson,
{\em Orthogonal functions from Gram determinants},
SIAM J. Math. Anal. 22 (1991), 1147--1155.
%
\bibitem{Z}
A. S. Zhedanov,
``Hidden symmetry'' of Askey--Wilson polynomials,
\emph{Theoret. and Math. Phys.} \textbf{89} (1991), 1146--1157.
%
\end{thebibliography}
\end{document}